\input amstex\documentstyle {amsppt}  
\pagewidth{12.5 cm}\pageheight{19 cm}\magnification\magstep1
\topmatter
\title Parabolic character sheaves, II\endtitle
\author G. Lusztig\endauthor
\address Department of Mathematics, M.I.T., Cambridge, MA 02139\endaddress
\dedicatory{Dedicated to Boris Feigin on the occasion of his 50th birthday}
\enddedicatory
\thanks Supported in part by the National Science Foundation\endthanks
\endtopmatter   
\document 
\define\lt{\ltimes}
\define\opp{\bowtie}
\define\dsv{\dashv}

\define\Lie{\text{\rm Lie }}
\define\po{\text{\rm pos}}
\define\frl{\forall}

\define\si{\sim}

\define\sqc{\sqcup}

\define\qua{\quad}

\define\hG{\hat G}
\define\hX{\hat X}

\define\tcl{\ti\cl}
\define\tSi{\ti\Si}
\define\tSS{\ti\SS}
\define\bA{\bar A}

\define\bG{\bar G}

\define\bZ{\bar Z}
\define\bvt{\bar\vt}

\define\lb{\linebreak}

\define\op{\oplus}

\define\em{\emptyset}

\define\iy{\infty}
\define\m{\mapsto}
\define\do{\dots}

\define\bsl{\backslash}

\define\lra{\leftrightarrow}

\define\sub{\subset}

\define\T{\times}
\define\ti{\tilde}
\define\nl{\newline}
\redefine\i{^{-1}}
\define\fra{\frac}
\define\un{\underline}

\define\bbq{\bar{\QQ}_l}

\define\Ad{\text{\rm Ad}}

\define\End{\text{\rm End}}

\define\supp{\text{\rm supp}}

\define\a{\alpha}
\redefine\b{\beta}

\define\g{\gamma}
\redefine\d{\delta}

\define\io{\iota}
\redefine\o{\omega}
\define\p{\pi}

\define\r{\rho}
\define\s{\sigma}
\redefine\t{\tau}
\define\th{\theta}
\define\k{\kappa}

\define\vt{\vartheta}

\redefine\G{\Gamma}
\redefine\D{\Delta}

\define\Si{\Sigma}

\define\Ph{\Phi}

\define\kk{\bold k}

\redefine\ss{\bold s}
\redefine\tt{\bold t}

\define\zz{\bold z}
\redefine\xx{\bold x}

\define\FF{\bold F}

\define\KK{\bold K}

\define\QQ{\bold Q}

\define\SS{\bold S}

\define\cb{\Cal B}
\define\cc{\Cal C}
\define\cd{\Cal D}

\define\ch{\Cal H}

\define\ck{\Cal K}
\define\cl{\Cal L}

\define\cp{\Cal P}

\define\cs{\Cal S}
\define\ct{\Cal T}

\define\cv{\Cal V}

\define\cz{\Cal Z}
\define\cx{\Cal X}
\define\cy{\Cal Y}

\define\tb{\ti b}

\define\tg{\ti g}

\define\tu{\ti u}

\define\tz{\ti z}

\define\tB{\ti B}

\define\tJ{\ti J}

\define\tP{\ti P}  
\define\tQ{\ti Q}

\define\tX{\ti X}
\define\tY{\ti Y}

\define\sh{\sharp}

\define\sps{\supset}
\define\BE{B}
\define\BBD{BBD}
\define\DP{DP}
\define\CS{L3}
\define\PCS{L9}
\head Introduction\endhead
Let $G$ be a connected reductive algebraic group over an algebraically closed field
$\kk$. Let $Z$ be the algebraic variety consisting of all triples $(P,P',U_{P'}gU_P)$ 
where $P,P'$ run through some fixed conjugacy classes of parabolics in $G$ and $g$ is 
an element of $G$ that conjugates $P$ to a parabolic in a fixed "good" relative 
position $y$ with $P'$ (here $U_P,U_{P'}$ are the unipotent radicals of $P,P'$). The 
varieties $Z$ include more or less as a special case the boundary pieces of the De 
Concini-Procesi completion $\bG$ of $G$ (assumed to be adjoint). They also include as a
special case the varieties studied in the first part of this series \cite{\PCS} (where 
$y=1$ that is, $gPg\i=P'$). In this special case a theory of "character sheaves" on $Z$
was developed in \cite{\PCS}. In the present paper we extend the theory of character 
sheaves to a general $Z$. 

We now review the content of this paper in more detail. (The numbering of sections
continues that of \cite{\PCS}; we also follow the notation of \cite{\PCS} .)

In Section 8 we introduce a partition of $Z$ similar to that in \cite{\PCS}; as in
\cite{\PCS}, it is based on the combinatorics in Section 2. But whereas in \cite{\PCS}
the combinatorics needed is covered by the results in \cite{\BE}, for the present paper
we actually need the slight generalization of \cite{\BE} given in Section 2. Now, it is
not obvious that the partition of $Z$ defined in Section 8 reduces for $y=1$ to that in
Section 3; this needs an argument that is given in Section 9. In Section 10 we consider
the example where $G$ is a general linear group. In Section 11 we define the "parabolic
character sheaves" on $Z$. As in the case $y=1$ (Section 4), we give two definitions 
for these; one uses the partition in Section 8 and allows us to enumerate the parabolic
character sheaves; the other one imitates the definition of character sheaves in
\cite{\CS}. (The two definitions are equivalent by 11.15 and 11.18.) The theory of
character sheaves in Section 11 generalizes that in Section 4 (this is seen from the
second definition). A consequence of the coincidence of the two definitions of
parabolic character sheaves on $Z$ is that a statement like 0.1(a) (concerning 
characteristic functions over a finite field) continues to hold in the generality of 
this paper. In Section 12 we define the notion of character sheaf on the completion
$\bG$. We again have two definitions; one is based on the partition in Section 8, and
the second one is reminiscent of the definition of character sheaves in \cite{\CS}. We
expect that these two definitions coincide, but we cannot prove this; if this was true,
we would have an analogue of 0.1(a) for $\bG$ over a finite field. Our results also 
provide a finite partition of $\bG$ into $G$-stable pieces (a refinement of the usual 
partition of $\bG$ into $G\T G$-orbits) which allows us to give an explicit description
of the set of $G$-conjugacy classes in $\bG$ (see 12.3(a)).

\head Contents\endhead
8. The variety $Z_{J,y,\d}$ and its partition.

9. Comparison of two partitions.

10. Example.

11. Parabolic character sheaves on $Z_{J,y,\d}$.

12. Completion. 

\head 8. The variety $Z_{J,y,\d}$ and its partition\endhead
\subhead 8.1 \endsubhead
We preserve the setup of 3.1. (Thus, $\hG$ is a possibly disconnected reductive 
algebraic group over $\kk$ with identity component $G$ and $G^1$ is a fixed connected 
component of $\hG$. Also $W,I$ is the Weyl group of $G$ and $\d:W@>\si>>W$.) Let 
$P\in\cp_J,Q\in\cp_K,u=\po(P,Q)$. We have 
$$\dim((U_P\cap U_Q)\bsl U_P)=l(u)+\nu_J-\nu_{J\cap\Ad(u)K}.\tag a$$
Here $\nu_J$ is the number of reflections in $W_J$. 

\subhead 8.2\endsubhead
Let $P,P'$ be two parabolics of $G$. The following hold. 

(a) $P^{P'},P'{}^P$ are in good position and $\po(P,P')=\po(P^{P'},P'{}^P)$;

(b) if $B\in\cb,B\sub P^{P'}$ then for some $B'\in\cb,B'\sub P'$ we have
$\po(B,B')=\po(P,P')$.
\nl
To prove (b) we may replace $P,P'$ by $P^{P'},P'{}^P$. It suffices to prove: if $P,P'$ 
are in good position and $B\sub P$, then for some $B'\sub P'$ we have 
$\po(B,B')=\po(P,P')$.

\subhead 8.3\endsubhead
If $P,Q$ are parabolics in good position, we have a bijection
$$\{\text{parabolics contained in }P\}@>\si>>\{\text{parabolics contained in }Q\}$$
given by $P'\lra Q'$, $Q'=Q^{P'},P'=P^{Q'}$. (Then $P',Q'$ are in good position and
$\po(P',Q')=\po(P,Q)$.)

\proclaim{Lemma 8.4}Let $P,Q,R$ be parabolics with a common Levi $L$. Then \lb
$\po(P,Q)\po(Q,R)=\po(P,R)$.
\endproclaim
Let $\b$ be a Borel of $L$. Then $B=U_P\b,B'=U_Q\b,B''=U_R\b$ are Borels of $P,Q,R$ 
respectively and we have
$$\po(P,Q)=\po(B,B'),\po(Q,R)=\po(B',B''),\po(P,R)=\po(B,B'').$$
It suffices to show that $\po(B,B')\po(B',B'')=\po(B,B'')$. This holds since 
$B,B',B''$ contain $\b$ hence have a common maximal torus. 

\proclaim{Lemma 8.5}Let $P,Q,R$ be parabolics with a common Levi $L$; let $P',Q',R'$ be
parabolics with a common Levi $L'$. Assume that $\po(P,Q)=\po(P',Q')$, 
$\po(Q,R)=\po(Q',R')$. Assume also that $P,P'$ have the same type; $Q,Q'$ have the same
type; $R,R'$ have the same type. Then there exists $x\in G$ that conjugates $P,Q,R$ to 
$P',Q',R'$.
\endproclaim
Clearly, we can assume that $P=P',Q=Q',L=L'$. Then we use the following fact: if 
$Q,R,R'$ are parabolics with a common Levi $L$ such that $\po(Q,R)=\po(Q,R')$ then
$R=R'$. (This can be reduced to the case where $Q,R,R'$ are Borels, which is clear.)

\proclaim{Lemma 8.6} Let $Q,Q'$ be parabolics with a common Levi. Then 
$Q\cap U_{Q'}=U_Q\cap U_{Q'}$.
\endproclaim
It suffices to show that $\Lie Q\cap\Lie U_{Q'}=\Lie U_Q\cap\Lie U_{Q'}$. Let $L$ be a
common Levi of $Q,Q'$. Consider the weight decomposition $\Lie G=\op_\a\Lie G_\a$ with 
respect to the connected centre of $L$. Then $\Lie G_0=\Lie L$ and 
$\Lie U_Q,\Lie U_{Q'}$ are direct sums of various $\Lie G_\a$ with $\a\ne 0$. Let 
$x\in\Lie Q\cap\Lie U_{Q'}$. We have $x=x_0+x'$ with $x_0\in\Lie L'$ and 
$x'\in\Lie U_Q$ is in $\op_{\a\ne 0}\Lie G_\a$. Since $x\in\Lie U_{Q'}$, we have 
$x\in\op_{\a\ne 0}\Lie G_\a$. Hence $x_0=0$ and $x\in\Lie U_Q$. Thus, 
$\Lie Q\cap\Lie U_{Q'}\sub\Lie U_Q\cap\Lie U_{Q'}$. The reverse inclusion is obvious.

\proclaim{Lemma 8.7} Let $P,P',Q$ be parabolics such that $\po(P',P)=z$. Assume that 
$\po(P',Q)=y$, that $P',Q$ are in good position and that $Q$ contains a Levi of 
$P\cap P'$. Let $Q'=Q^{(P'{}^P)}$. Then

(a) $\po(P^{P'},P'{}^P)=z\i$.

(b) $\po(P'{}^P,Q')=y$.

(c) $\po(P^{P'},Q')=z\i y$ and $P^{P'},Q'$ are in good position.
\endproclaim
(a) follows from 8.2; (b) follows from 8.3. We prove (c). Let $L_0$ be a common Levi of
$P^{P'},P'{}^P$ that is contained in $Q$. Then $P^{P'},P'{}^P,Q'$ have a common Levi 
$L_0$. By Lemma 8.4 we have $\po(P^{P'},Q')=\po(P^{P'},P'{}^P)\po(P'{}^P,Q')=z\i y$ as
required. 

\subhead 8.8\endsubhead
Let $J,J'\sub I$ and $y\in W$ be such that $\Ad(y)\d(J)=J',y\in{}^{J'}W^{\d(J)}$. For
$P,P'$ in $\cp_J,\cp_{J'}$ let 

$A_y(P,P')=\{g\in G^1;\po(P',{}^gP)=y\}$,
 
$A'_y(P,P')=\{g\in A_y(P,P');{}^gP\text{ contains some Levi of }P\cap P'\}$.

\proclaim{Lemma 8.9}(a) $A_y(P,P')$ is a single $P',P$ double coset and also a single 
$P',U_P$ double coset.

(b) $A_y(P,P')=U_{P'}A'_y(P,P')$.
 
(c) $A'_y(P,P')$ is a single $(P\cap P'),P$ double coset and also a single
$U_{P\cap P'},P$ double coset.
\endproclaim
We prove (a). We can find $Q\in\cp_{\d(J)}$ such that $\po(P',Q)=y$; moreover, we can
find $g\in G^1$ such that ${}^gP=Q$. Thus $A_y(P,P')\ne\em$. Let $g,g'\in A_y(P,P')$. 
Clearly, $g'=xgp$ with $x\in P',p\in P$. Now ${}^gP,P'$ are in good position; let $L$ 
be a common Levi of them. Since ${}^{g\i}L$ is a Levi of $P$, we have 
$p\in{}^{g\i}LU_P$. Thus $p=g\i lgu$ with $l\in L,u\in U_P$ and 
$g'=xgp=xlgu\in P'gU_P$.

We prove (b). Let $g\in A_y(P,P')$. Let $L_0$ be a Levi of $P\cap P'$. Then $L_0$ is
contained in a Levi $L_1$ of $P'$. Let $L_2$ be a common Levi of $P',{}^gP$. Then 
$L_1={}^uL_2$ for some $u\in U_{P'}$. We have
${}^{ug}P={}^u({}^gP)\sps{}^uL_2=L_1\sps L_0$ hence $ug\in A'_y(P,P')$.

We prove (c). Let $g,g'\in A'_y(P,P')$. Then 

$P^{P'},P'{}^P,({}^gP)^{(P'{}^P)}$ have a common Levi;

$P^{P'},P'{}^P,({}^{g'}P)^{(P'{}^P)}$ have a common Levi;

$\po(P'{}^P,({}^gP)^{(P'{}^P)})=\po(P'{}^P,({}^{g'}P)^{(P'{}^P)})=y$.
\nl
By Lemma 8.5, there exists $x\in G$ which conjugates 

$P^{P'},P'{}^P,({}^gP)^{(P'{}^P)}$ to $P^{P'},P'{}^P,({}^{g'}P)^{(P'{}^P)}$. 
\nl
Then $x\in P^{P'}\cap P'{}^P=P\cap P'$ and $x$ conjugates ${}^gP$ to ${}^{g'}P$, since 
${}^gP$ to ${}^{g'}P$ are parabolics of type $\d(J)$ containing 
$({}^gP)^{(P'{}^P)},({}^{g'}P)^{(P'{}^P)}$. Hence $xg\in g'P$ that is $g'\in xgP$. Let
$M$ be a Levi of $P\cap P'$ with $M\sub{}^gP$. We can write $x=vm$ with
$v\in U_{P\cap P'},m\in M$. Then $g\i mg\in P$, $xgP=vmgP=vgg\i mgP=vgP$. The lemma is 
proved.

\subhead 8.10\endsubhead
Let $P\in\cp_J,P'\in\cp_{J'}$ be such that $\po(P',P)=z$. Let 
$$J_1=J\cap \d\i\Ad(y\i z)J,J'_1=J\cap\Ad(z\i y)\d(J).$$
Then $\Ad(z\i y)\d(J_1)=J'_1$. Let $g\in A'_y(P,P')$. We set 
$$P_1=g\i({}^gP)^{(P'{}^P)}g,\qua P'_1=P^{P'},\qua (P_1,P'_1)=\a(P,P',g).$$
We have $P_1\sub P,P'_1\sub P$ and, by Lemma 8.7, $P'_1,{}^gP_1$ are in good position,
$$P_1\in\cp_{J_1},P'_1\in\cp_{J'_1},\qua\po(P'_1,{}^gP_1)=z\i y.$$
(We have also $\po(P'_1,P'{}^P)=z\i,\po(P'{}^P,{}^gP_1)=y$). Thus, 
$z\i y\in{}^{J'_1}W^{\d(J_1)}$ and $g\in A_{z\i y}(P_1,P'_1)$.

\proclaim{Lemma 8.11} Let $g,g'\in A'_y(P,P'),u'\in U_{P'},u\in U_P$ with $g'=u'gu$.
Then 

(a) $g'=u'_1gu_1$ with $u'_1\in U_{P'_1},u_1\in U_{P_1}$;

(b) we have $\a(P,P',g')=(P_1,P'_1)$.
\endproclaim
We prove (a). Since $P_1\sub P$ we have $U_P\sub U_{P_1}$. Hence we may assume that
$u=1$. By Lemma 8.9(c) we have $g'=vgp$ with $v\in U_{P\cap P'},p\in P$. Thus 
$g'=vgp=u'g$. Hence $v\i u'=gpg\i$. Now $v\i u'\in U_{P\cap P'}U_{P'}\sub U_{P'{}^P}$. 
Thus $v\i u'\in U_{P'{}^P}\cap{}^gP$. We have 
$$U_{P'{}^P}\cap{}^gP=U_{P'{}^P}\cap{}^gP_1=U_{P'{}^P}\cap U_{{}^gP_1}.$$
(In general, if $R,S$ are parabolics in good position and $R'\sub R,S'=S^{R'}$, then 
$U_{R'}\cap S=U_{R'}\cap S'$. Indeed, $S'=(R'\cap S)U_S$ hence $U_{R'}\cap S\sub S'$. 
Since $R',S'$ are in good position, we have $U_{R'}\cap S'=U_{R'}\cap U_{S'}$ by Lemma 
8.6.) Thus, $gpg\i=v\i u'\in U_{{}^gP_1}$ so that $p\in U_{P_1}$. Now 
$v\in U_{P\cap P'}\sub U_{P^{P'}}=U_{P'_1}$. We see that $g'=vgp$ with 
$v\in U_{P'_1},p\in U_{P_1}$. This proves (a).

We prove (b). We must show that
$$u\i g\i u'{}\i({}^{u'gu}P)^{(P'{}^P)}u'gu=g\i({}^gP)^{(P'{}^P)}g$$
or that $u\i g\i({}^gP)^{(P'{}^P)}gu=g\i({}^gP)^{(P'{}^P)}g$. This holds since 
$g\i({}^gP)^{(P'{}^P)}g$ is a parabolic subgroup of $P$ hence it contains $u$. 

\proclaim{Lemma 8.12} Let $g,g'\in A'_y(P,P')$. Assume that 
$\a(P,P',g)=\a(P,P',g')=(P_1,P'_1)$ and that $g'\in U_{P'_1}gU_{P_1}$. Then there exist
$x\in U_P\cap P',w'\in U_{P'},w\in U_P$ such that $g'=w'xgw$.
\endproclaim
By Lemma 8.9(c) we have $g'=vgp$ with $v\in U_{P\cap P'},p\in P$. By assumption,
$p\i g\i v\i({}^{vgp}P)^{(P'{}^P)}vgp=g\i({}^gP)^{(P'{}^P)}g$ that is 
${}^{p\i}P_1=P_1$, or $p\in P_1$. Also $g'=u'gu$ with $u'\in U_{P'_1},u\in U_{P_1}$. 
Thus, $g'=vgp=u'gu$. Setting $\p=pu\i\in P_1$ we have $vg\p=u'g$ and $v\i u'={}^g\p$. 
Now $v\i u'\in U_{P\cap P'}U_{P^{P'}}=U_{P^{P'}}$. Thus 
$v\i u'\in U_{P^{P'}}\cap{}^gP_1$. Since $P^{P'},{}^gP_1$ are in good position, we have
$v\i u'\in U_{P^{P'}}\cap U_{{}^gP_1}$. Thus ${}^g\p\in U_{{}^gP_1}$ and
$\p\in U_{P_1}$. Since $u\in U_{P_1}$, we have $p\in U_{P_1}$. Thus 
$g'\in U_{P\cap P'}gU_{P_1}=U_{P\cap P'}U_{{}^gP_1}g$. Now ${}^gP_1=({}^gP)^{(P'{}^P)}$
hence $U_{{}^gP_1}=({}^gP\cap U_{P'{}^P})U_{{}^gP}$ so that
$$\align&g'\in U_{P\cap P'}({}^gP\cap U_{P'{}^P})U_{{}^gP}g=
U_{P\cap P'}({}^gP\cap U_{P'{}^P})gU_P\sub U_{P'{}^P}gU_P\\&
=U_{P'}(U_P\cap P')gU_P.\endalign$$
Thus, $g'=w'xgw$ with $w'\in U_{P'},x\in P'\cap U_P,w\in U_P$, as desired.

\subhead 8.13\endsubhead
We fix $z\in{}^{J'}W^J$. Let $J_1=J\cap\d\i\Ad(y\i z)J,J'_1=J\cap\Ad(z\i y)\d(J)$,
so that $\Ad(z\i y)\d(J_1)=J'_1$. Let $Q,Q'$ in $\cp_{J_1},\cp_{J'_1}$ be such that 
$\po(Q',Q)\in W_J$. Let $\g_1$ be a $U_{Q'},U_Q$ double coset in $A_{z\i y}(Q,Q')$. Let
$F'$ be the set of all $(P,P',g)$ with $P\in\cp_J,P'\in\cp_{J'},\po(P',P)=z$,
$g\in A'_y(P,P')$ such that $\a(P,P',g)=(Q,Q')$ and $g\in\g_1$. (Since $Q'\sub P$, $P$ 
is uniquely determined.) Let $F$ be the set of all $(P,P',\g)$ with
$P\in\cp_J,P'\in\cp_{J'},\po(P',P)=z$, $\g\in U_{P'}\bsl A_y(P,P')/U_P$ such that for
some/any $g\in\g\cap A'_y(P,P')$ we have $\a(P,P',g)=(Q,Q')$ and $g\in\g_1$. (The 
equivalence of "some/any" follows from Lemma 8.11.) Again, $P$ is uniquely determined. 
The map $F'@>>>F,(P,P',g)\m(P,P',U_{P'}gU_P)$ is surjective by Lemma 8.9(b).

\proclaim{Lemma 8.14} $F'$ is non-empty. 
\endproclaim
Let $g\in\g_1$. Then $\po(Q',{}^gQ)=z\i y$. Define $P\in\cp_J$ by $Q\sub P$. Since
$\po(Q',Q)\in W_J$, we have also $Q'\sub P$. Let $\tP,P'$ in $\cp_J,\cp_{J'}$ with
$\po(P',\tP)=z$. By Lemma 8.9(c), we can find $\tg\in A'_y(\tP,P')$. Now 
$\a(\tP,P',\tg)=(\tQ,\tQ')$ with $\tQ\in\cp_{J_1},\tQ'\in\cp_{J'_1}$ contained in $\tP$
and $\po(\tQ',{}^{\tg}\tQ)=z\i y$. We have $Q'={}^h\tQ'$ for some $h\in G$. Replacing 
$\tP,P',\tg$ by ${}^h\tP,{}^hP',h\tg h\i$, we may assume that $\tQ'=Q'$. Then
$\tP\in\cp_J$ contains $Q'$ hence $\tP=P$. Now $\tQ,Q$ are contained in $P$ and have 
the same type hence $\tQ={}^pQ$ for some $p\in P$. We have 
${}^{\tg p}Q=({}^{\tg}P)^{(P'{}^P)}$ and $\tg p\in A'_y(P,P')$ hence 
$\a(P,P',\tg p)=(Q,Q')$. Since $\po(Q',{}^gQ)=\po(Q',{}^{\tg p}Q)$ and $Q',{}^gQ$ are
in good position, we have $\tg p=u'gq$ with $q\in Q,u'\in U_{Q'}$. Thus 
$u'gq\in A'_y(P,P')$ and $\a(P,P',u'gq)=(Q,Q')$. We have $u'g\in A'_y(P,P')$ and 
$\a(P,P',u'g)=(Q,Q')$ (we use that $q\in Q$). Since $u'g\in\g_1$, we see that 
$F'\ne\em$.

\proclaim{Lemma 8.15} Let $(P,P',g)\in F'$. Then $(u,v)\m(P,{}^uP',U_{{}^uP'}uvgU_P)$ 
is a well defined, surjective map $\k:U_P\T(U_P\cap P')@>>>F$.
\endproclaim
Let $u\in U_P,v\in U_P\cap P'$. Clearly, $(P,{}^uP',uvg)\in F'$. Hence $\k$ is well
defined. We show that $\k$ is surjective. Let $(P,\tP',\g)\in F$. We can find 
$\tg\in\g$ such that $(P,\tP',\tg)\in F'$. We have $g\in\g_1,\tg\in\g_1$ hence 
$\tg\in U_{Q'}gU_Q$. Since $\po(P',P)=\po(\tP',P)$, we have $\tP'={}^pP'$ for some
$p\in P$. Since $P^{P'}=P^{\tP'}=P^{({}^pP')}=pP^{P'}p\i$, we have $p\in P^{P'}$. Thus,
$p=u\p$ with $\p\in P\cap P',u\in U_P$. Since ${}^pP'={}^uP'$, we may assume that
$p=u\in U_P$ so that $\tP'={}^uP'$. Applying Lemma 8.12 to 
$(P,{}^uP',ug),(P,{}^uP',\tg)\in F'$ (instead of $(P,P',g),(P,P',g')$) we see that
$\tg=w'xugw$ for some $x\in U_P\cap{}^uP',w'\in U_{{}^uP'},w\in U_P$. Let 
$v=u\i xu\in U_P\cap P'$. Then $\tg=w'uvgw$. Thus, $(P,\tP',\g)=\k(u,v)$. The lemma is 
proved.

\proclaim{Lemma 8.16} In the setup of Lemma 8.15, the following two conditions for
$(u,v),(u,v')$ in $U_P\T(U_P\cap P')$ are equivalent:

(i) $\k(u,v)=\k(u',v')$;

(ii) $u'=uf,v'=f\i dv$ for some $f\in U_P\cap P',d\in U_P\cap U_{P'}$.
\endproclaim
Assume that (i) holds. We have ${}^uP'={}^{u'}P',uvg\in U_{{}^uP'}u'v'gU_P$. Thus
$u'=uf$ with $f\in P'$ (hence $f\in U_P\cap P'$) and $uvg\in uU_{P'}fv'gU_P$, that is,
$v\in U_{P'}fv'U_{{}^gP}$, so that $v\in fv'U_{P'}U_{{}^gP}$. We show that
$$P'\cap U_{P'}U_{{}^gP}=U_{P'}.\tag a$$
Assume that $x\in P'\cap U_{P'}U_{{}^gP}$. We must show that $x\in U_{P'}$. We have
$u_1x\in U_{{}^gP}$ with $u_1\in U_{P'}$. Let $x'=u_1x\in P'$. Then 
$x'\in P'\cap U_{{}^gP}=U_{P'}\cap U_{{}^gP}$ (by Lemma 8.6, which is applicable since
$P',{}^gP$ have a common Levi). Thus, $x'\in U_{P'}$ hence $x\in U_{P'}$, as required.

Applying (a) to $v'{}\i f\i v\in P'\cap U_{P'}U_{{}^gP}$, we see that
$v'{}\i f\i v\in U_{P'}$, so that $v'=f\i dv$ for some $d\in U_{P'}$. We have 
$d=fv'v\i\in U_P\cap P'$. Hence $d\in U_P\cap U_{P'}$. Thus, (ii) holds. The converse 
is immediate. The lemma is proved.

\subhead 8.17\endsubhead
We consider a new group structure $(d,f)\bullet(d',f')=(f'df'{}\i d',f'f)$ on 
$(U_P\cap U_{P'})\T(U_P\cap P')$ and a new group structure 
$(u,v)\bullet(u',v')=(u'u,vv')$ on $U_P\T(U_P\cap P')$. Then 
$$\th:(U_P\cap U_{P'})\T(U_P\cap P')@>>>U_P\T(U_P\cap P'),\qua (d,f)\m(f,f\i d),$$
is an (injective) group homomorphism for these new group structures. We can reformulate
condition (ii) in 8.16 as follows:
$$(u',v')=\th(x)\bullet(u,v)\text{ for some }x\in(U_P\cap U_{P'})\T(U_P\cap P').$$
We see that $\k$ defines a bijection 
$$\th(((U_P\cap U_{P'})\T(U_P\cap P'))\bsl(U_P\T(U_P\cap P')))@>\si>>F.$$
One can check that this is in fact an isomorphism of algebraic varieties. Since 
$U_P\T(U_P\cap P')$ is a connected unipotent group and 
$\th((U_P\cap U_{P'})\T(U_P\cap P'))$ is a connected closed subgroup of it, we see that

(a) $F$ is isomorphic to an affine space of dimension $\dim(U_P/(U_P\cap U_{P'}))$.

\subhead 8.18\endsubhead
Let $J,J'\sub I$ and $y\in{}^{J'}W^{\d(J)}$ be such that $\Ad(y)\d(J)=J'$. Let
$$Z_{J,y,\d}=\{(P,P',\g);P\in\cp_J,P'\in\cp_{J'},\g\in U_{P'}\bsl A_y(P,P')/U_P\}.$$
To any $(P,P',\g)\in Z_{J,y,\d}$ we associate a sequence $(J_n,J'_n,u_n)_{n\ge 0}$ with
$J_n,J'_n\sub I$, $u_n\in W$, a sequence $(y_n)_{n\ge 0}$ with 
$y_n\in{}^{J'_n}W^{\d(J_n)},\Ad(y_n)\d(J_n)=J'_n$ and a sequence 
$(P_n,P'_n,\g_n)_{n\ge 0}$ with 
$P_n\in\cp_{J_n},P'_n\in\cp_{J'_n},\g_n\sub A_{y_n}(P_n,P'_n)$. We set 
$$P_0=P,P'_0=P',\g_0=\g,J_0=J,J'_0=J',u_0=\po(P'_0,P_0),y_0=y.$$  
Assume that $n\ge 1$, that $P_m,P'_m,\g_m,J_m,J'_m,u_m,y_m$ are already defined for 
$m<n$ and that $u_m=\po(P'_m,P_m),P_m\in\cp_{J_m},P'_m\in\cp_{J'_m}$ for $m<n$. Let 
$$J_n=J_{n-1}\cap\d\i\Ad(y_{n-1}\i u_{n-1})J_{n-1},
J'_n=J_{n-1}\cap\Ad(u_{n-1}\i y_{n-1})\d(J_{n-1}),$$
$$P_n=g_{n-1}\i({}^{g_{n-1}}P_{n-1})^{(P'_{n-1}{}^{P_{n-1}})}g_{n-1}\in\cp_{J_n},
P'_n=P_{n-1}^{P'_{n-1}}\in\cp_{J'_n}$$
where 
$$g_{n-1}\in\g_{n-1}\cap A'_{y_{n-1}}(P_{n-1},P'_{n-1}),$$
$$u_n=\po(P'_n,P_n),y_n=u_{n-1}\i y_{n-1},\g_n=U_{P'_n}g_{n-1}U_{P_n}.$$
This completes the inductive definition; the definition makes sense (it is independent 
of choices) by 8.9-8.11. We have $(J_n,J'_n,u_n)_{n\ge 0}\in S(J,\Ad(y)\d)$ (see 2.3).
We write $(J_n,J'_n,u_n)_{n\ge 0}=\b(P,P',\g)$. For $\ss\in S(J,\Ad(y)\d)$ let 
$$Z_{J,y,\d}^\ss=\{(P,P',\g)\in Z_{J,y,\d};\b(P,P',\g)=\ss\}.$$
Clearly, $(Z_{J,y,\d}^\ss)_{\ss\in S(J,\Ad(y)\d)}$ is a partition of $Z_{J,y,\d}$ into 
locally closed subvarieties. The $G$-action on $Z_{J,y,\d}$ given by
$g:(P,P',\g)\m({}^gP,{}^gP',g\g g\i)$ preserves each of the pieces $Z_{J,y,\d}^\ss$.
Now $(P,P',\g)\m(P_1,P'_1,\g_1)$ is a morphism
$f:Z_{J,y,\d}^\ss@>>>Z_{J_1,y_1,\d}^{\ss_1}$ where for 
$\ss=(J_n,J'_n,u_n)_{n\ge 0}\in S(J,\Ad(y)\d)$ we set 
$\ss^1=(J_n,J'_n,u_n)_{n\ge 1}\in S(J_1,\Ad(y_1)\d)$.

\proclaim{Lemma 8.19}(a) The morphism $\vt:Z_{J,y,\d}^\ss@>>>Z_{J_1,y_1,\d}^{\ss_1}$ is
a locally trivial fibration with fibres isomorphic to an affine space of dimension 
$l(u_0)+\nu_J-\nu_{J_1}$.

(b) Let $\bvt$ be the map from the set of $G$-orbits on $Z_{J,y,\d}^\ss$ to the set 
of $G$-orbits on $Z_{J_1,y_1,\d}^{\ss_1}$ induced by $\vt$. Then $\bvt$ is a bijection.
\endproclaim
We prove (a). Let $(P,P',\g)\in Z_{J,y,\d}^\ss$. From 8.17 and 8.1(a) we see that each 
fibre of $\vt$ is an affine space of dimension
$$\dim(U_P/(U_P\cap U_{P'}))=l(u)+\nu_J-\nu_{J\cap\Ad(u_0\i)J'}=l(u)+\nu_J-\nu_{J'_1}
=l(u)+\nu_J-\nu_{J_1}.$$
The verification of local triviality is omitted.

We prove (b). From the fact that $\vt$ is surjective (see (a)) and $G$-equivariant, it
follows that $\bvt$ is well defined and surjective. We show that $\bvt$ is injective. 
Let $(P,P',\g),(\tP,\tP',\ti\g)$ be two triples in $Z_{J,y,\d}^\ss$ whose images under 
$\vt$ are in the same $G$-orbit; we must show that these two triples are in the same 
$G$-orbit. Since $\vt$ is $G$-equivariant, we may assume that 
$\vt(P,P',\g)=\vt(\tP,\tP',\ti\g)=(Q,Q',\g_1)\in Z_{J_1,y_1,\d}^{\ss_1}$. Define $F$ in
terms of $(Q,Q',\g_1)$ as in 8.13. Then $(P,P',\g)\in F,(\tP,\tP',\ti\g)\in F$. Since 
$P,\tP$ are parabolics of the same type containing $Q'$ we have $P=\tP$. Let 
$g\in\g\cap A'_y(P,P')$. By Lemma 8.15, there exist $u\in U_P,v\in U_P\cap P'$ such
that $\tP'={}^uP'$, $\ti\g=U_{{}^uP'}uvgU_P$. We have also $\tP={}^{uv}P$ (since 
$uv\in P$), $\tP'={}^{uv}P'$ (since $v\in P'$), 
$$\ti\g=uU_{P'}vgU_P=uvU_{P'}gU_P=uvU_{P'}gU_Pv\i u\i=uv\g v\i u\i$$
(since $v$ normalizes $U_{P'}$ and $uv\in U_P$). Thus, $(\tP,\tP',\ti\g)$ is obtained 
by the action of $uv\in G$ on $(P,P',\g)$, hence $(\tP,\tP',\ti\g)$ is in the $G$-orbit
of $(P,P',\g)$. The lemma is proved.

\proclaim{Lemma 8.20} Let $\ss=(J_n,J'_n,u_n)_{n\ge 0}\in S(J,\Ad(y)\d)$. Then
$Z_{J,y,\d}^\ss$ is an iterated affine space bundle (with fibre dimension 
$l(w)+\nu_{J}-\nu_{J_m}$, $w=u_0u_1\do u_m$, $m\gg 0$) over a fibre bundle over
$\cp_{J_m}$ with fibres isomorphic to $P/U_P$ with $P\in\cp_{J_m}$, $m\gg 0$. In
particular, $Z_{J,y,\d}^\ss\ne\em$.
\endproclaim
Assume first that $\ss$ is such that $J_n=J'_n=J$ and $u_n=1$ for all $n\ge 0$. (Then 
$\Ad(y)\d(J)=J,y\in{}^JW^{\d(J)}$.) In this case, $Z_{J,y,\d}^\ss$ is the set of all 
$(P,P',\g)$ with $P=P'\in\cp_J$, $\g\in U_P\bsl A_y(P,P)/U_P$. (The associated 
sequence $(P_n,P'_n,\g_n)$ is in this case $P_n=P'_n=P,\g_n=\g$.) Thus, 
$Z_{J,y,\d}^\ss$ is a locally trivial fibration over $\cp_J$ with fibres isomorphic to 
$P/U_P$ for $P\in\cp_J$ and the lemma holds.

We now consider a general $\ss$. For any $r\ge 0$ let 
$$\ss_r=(J_n,J'_n,u_n)_{n\ge r}\in S(J_r,\Ad(y_r)\d)$$ 
($y_r$ as in 8.18). By 8.19(a) we have a sequence of affine space bundles
$$Z_{J,y,\d}^\ss@>>>Z_{J_1,y_1,\d}^{\ss_1}@>>>Z_{J_2,y_2,\d}^{\ss_2}@>>>\do\tag a$$
where for $r\gg 0$, $Z_{J_r,y_r,\d}^{\ss_r}$ is as in the first part of the proof. By
8.19(a), the sum of dimensions of fibres of the maps in this sequence is
$$\sum_{n\ge 0}(l(u_n)+\nu_{J_n}-\nu_{J_{n+1}})=\sum_{n\ge 0}l(u_n)+\nu_{J_0}-\nu_{J_m}
                                               =l(w)+\nu_J-\nu_{J_m}$$
where $m\gg 0$. The lemma follows.

{\it Remark.} If $\hG$ is defined over the finite field $\FF_p$ and $\kk$ is the 
algebraic closure of $\FF_p$ then the number $N$ of rational points of $Z_{J,J',y}^\ss$
over a sufficiently large finite subfield $\FF_q$ of $\kk$ equals
$$\sh(G(\FF_q))q^{l(w)+\nu_J-\nu_I}.$$
Indeed, for $m\gg 0$ we have
$N=\sh(G(\FF_q)/P_m(\FF_q))\fra{\sh(P_m(\FF_q))}{q^{\nu_I-\nu_{J_m}}}
q^{l(w)+\nu_J-\nu_{J_m}}$. 
Note that $l(w)+\nu_J-\nu_I\ge 0$.

\subhead 8.21\endsubhead
In the setup of 8.20, the maps in 8.20(a) induce bijections on the sets of $G$-orbits 
(see 8.19(b)). Thus we obtain a canonical bijection between the set of $G$-orbits on 
$Z_{J,y,\d}^\ss$ and the set of $G$-orbits on $Z_{J_r,y_r,\d}^{\ss_r}$ with $r$ large 
enough so that $J_r=J'_r=J_{r+1}=J'_{r+1}=\do$, $u_r=u_{r+1}=\do=1$. This last set of 
orbits is canonically the set of $Q$-orbits on $U_Q\bsl A_{y_r}(Q,Q)/U_Q$ where 
$Q\in\cp_{J_r}$. The $Q$-action (conjugation) factors through $Q/U_Q$. Let $L^\ss$ be a
Levi of $Q$. Then 
$$C^\ss=\{g\in G^1;{}^gL^\ss=L^\ss,\po({}^gQ,Q)=y_r\}$$
is a connected component of $N_{\hG}(L^\ss)$. We have an obvious bijection 
$C^\ss@>\si>>U_Q\bsl A_{y_r}(Q,Q)/U_Q$ under which the action of $L^\ss$ on $C^\ss$ by 
conjugation corresponds to the action of $Q/U_Q$ on $U_Q\bsl A_{y_r}(Q,Q)/U_Q$ by 
conjugation. Thus we obtain a canonical bijection between the set of $G$-orbits on 
$Z_{J,y,\d}^\ss$ and the set of $L^\ss$-conjugacy classes in $C^\ss$ (a connected 
component of an algebraic group with identity component $L^\ss$). Putting together 
these bijections we obtain a bijection 
$$G\bsl Z_{J,y,\d}\lra\sqc_{\ss\in S(J,\Ad(y)\d)}L^\ss\bsl C^\ss\tag a$$
where $G\bsl Z_{J,y,\d}$ is the set of $G$-orbits on $Z_{J,y,\d}$ and $L^\ss\bsl C^\ss$
is the set of $L^\ss$-orbits on $C^\ss$ (for the conjugation action).

\head 9. Comparison of two partitions\endhead
\subhead 9.1\endsubhead
In the case where $y=1$, Sections 3 and 8 provide two methods to partition 
$Z_{J,y,\d}$. In this section we show that the resulting partitions of $Z_{J,y,\d}$ are
the same. Lemmas 9.2, 9.3 hold for any $y$, but in 9.4, 9.5 we assume that $y=1$.

\proclaim{Lemma 9.2} Let $(P,P',\g)\in Z_{J,y,\d}$. Let $n\ge 1$. Let $P'_1,P_n$ be as
in 8.18. We have $\po(P',P_n)=\po(P',P'_1)\po(P'_1,P_n)$.
\endproclaim
Let $z=\po(P',P'_1)$, $\tz=\po(P',P_n),x=\po(P'_1,P_n)$. We have $z=\po(P',P)$ hence 
$z\in{}^{J'}W^J$. We have also $\tz\in{}^{J'}W$. Since $P_n\sub P$ we have 
$\tz\in W_{J'}zW_J$. Using 2.1(c) with $x,x'$ replaced by $z,\tz$, we see that $\tz=zv$
with $v\in W_J$. Let $B,B'\in\cb$ be such that $B\sub P'$, $B'\sub P_n$,
$\po(B,B')=\tz$. Since $z\in W^J$, we have $l(zv)=l(z)+l(v)$. Hence there is a unique 
$B''\in\cb$ such that $\po(B,B'')=z$, $\po(B'',B')=v$. Since $B'\sub P$ and
$\po(B'',B)\in W_J$, we have $B''\sub P$. Since $B\sub P',B''\sub P$ and 
$\po(B,B'')=\po(P',P)=z$, we have $B''\sub P^{P'}=P'_1$. Since 
$B''\sub P'_1,B'\sub P_n$, we have $\po(B'',B')\ge x$, hence $v\ge x$. We can find
$B_1,B_2\in\cb$ such that $B_1\sub P'_1$, $B_2\sub P_n$, $\po(B_1,B_2)=x$. Since 
$\po(P',P)=z$ and $B_1\sub P^{P'}$, we can find $B_0\in\cb$ such that $B_0\sub P'$,
$\po(B_0,B_1)=z$. Since $z\in W^J,x\in W_J$, we have $\po(B_0,B_2)=zx$. We have 
$B_0\sub P'$, $B_2\sub P_n$ hence $\po(B_0,B_2)\ge\po(P',P_n)$, that is, 
$zx\ge\tz=zv$. Thus, we have $v\ge x$ and $zx\ge zv$. Since $z\in W^J$ and $x,v\in W_J$
we have $x=v$. Thus, $\tz=zx$. The lemma is proved.

\proclaim{Lemma 9.3} Let $(P,P',\g)\in Z_{J,y,\d}$. Let $P'_n,P_n,u_n,\g_n$ be as in
8.18. For any $n\ge 0$ we have $\po(P',P_n)=u_0u_1\do u_n$.
\endproclaim
We argue by induction on $n=0$. For $n=0$ the result is clear. Assume now that 
$n\ge 1$. We have $\po(P',P'_1)=\po(P',P^{P'})=\po(P',P)=u_0$. Using the induction 
hypothesis for $(P_1,P'_1,\g_1)$ and $n-1$ (instead of $(P,P',\g)$ and $n$) we see that
$\po(P'_1,P_n)=u_1\do u_n$. By Lemma 9.2 we have
$\po(P',P_n)=\po(P',P'_1)\po(P'_1,P_n)=u_0(u_1\do u_n)$. The lemma is proved.

\proclaim{Lemma 9.4} Let $J\sub I$. Let $(P,P',g)\in\cz_{J,\d}$. To $(P,P',g)$ we 
associate $P^n,P'{}^n,w_n$ as in 4.11. To $(P,P',U_{P'}gU_P)$ we associate 
$P_n,P'_n,\g_n,u_n,y_n$ as in 8.18 (with $y=1$). For any $n\ge 0$, the following hold:

($a_n$)  $P_n=P^n$;

($b_n$)  $P'_n=(P^{n-1})^{P'{}^{n-1}}$;

($c_n$) $g\in\g_n\cap A'_{y_n}(P_n,P'_n)$;

($d_n$) $(P'{}^n)^{(P'_n{}^{P^n})}=(P'{}^n)^{P^n}$;

($e_n$) $w_n=u_0u_1\do u_n$.
\endproclaim
We argue by induction on $n$. The result is obvious for $n=0$: we have $P_0=P^0=P$; 
${}^gP_0=P'$ contains a Levi of $P\cap P'$; we have $P'{}^{(P'{}^P)}=P'{}^P$. We have 
$w_0=u_0$. We now assume that $n\ge 1$ and that 
$(a_{n-1}),(b_{n-1}),(c_{n-1}),(d_{n-1}),(e_{n-1})$ hold.

We show that $(a_n)$ holds. Using $(c_{n-1}),(a_{n-1}),(d_{n-1})$, we have
$$\align&P_n=g\i({}^gP_{n-1})^{(P'_{n-1}{}^{P_{n-1}})}g=
g\i({}^gP^{n-1})^{(P'_{n-1}{}^{P^{n-1}})}g\\&
=g\i(P'{}^{n-1})^{(P'_{n-1}{}^{P^{n-1}})}g=g\i(P'{}^{n-1})^{P^{n-1}}g=P^n.\endalign$$
hence $(a_n)$ holds.

We show that $(b_n)$ holds, using $(a_{n-1}),(b_{n-1}),(e_{n-1}),(a_n)$. It suffices to
show that $Y'=P'_n$ satisfies the hypotheses of Lemma 3.2(a) with $P,P',a,J,\d(J)$ 
replaced by $P^{n-1},P'{}^{n-1},w_{n-1},J_{n-1},\d(J_{n-1})$. By definition we have 
$Y'\sub P_{n-1}=P^{n-1}$. The type of $Y'$ is 
$J_{n-1}\cap\Ad(u_0u_1\do u_{n-1})\i\d(J_{n-1})$ while that of $(P^{n-1})^{P'{}^{n-1}}$
is $J_{n-1}\cap\Ad(w_{n-1}\i)\d(J_{n-1})$. Thus, $Y'$ and $(P^{n-1})^{P'{}^{n-1}}$ have
the same type. Since ${}^gP_n={}^gP^n=P'{}^n=(P'{}^{n-1})^{P^{n-1}}$ and $P'_n,{}^gP_n$
are in good position, $\po(P'_n,{}^gP_n)=(u_0u_1\do u_{n-1})\i=w_{n-1}\i$, we see that
$(P'{}^{n-1})^{P^{n-1}},P'_n$ are in good position, 
$\po((P'{}^{n-1})^{P^{n-1}},P'_n)=w_{n-1}$. Thus, the hypotheses of Lemma 3.2(a) are
satisfied and $(b_n)$ holds.

We show that $(c_n)$ holds, using $(a_n),(b_n),(c_{n-1})$. Since \lb
$g\in\g_{n-1}\cap A'_{y_{n-1}}(P_{n-1},P'_{n-1})$ and 
$\g_{n-1}\cap A'_{y_{n-1}}(P_{n-1},P'_{n-1})\sub\g_n\cap A_{y_n}(P_n,P'_n)$ (see 8.10) 
we have $g\in\g_n\cap A_{y_n}(P_n,P'_n)$. It remains to show that ${}^gP_n$ contains a
Levi of $P_n\cap P'_n$ or equivalently, that $(P'{}^{n-1})^{P^{n-1}}$ contains a Levi 
of $(P^{n-1})^{P'{}^{n-1}}\cap g\i(P'{}^{n-1})^{P^{n-1}}g$. This follows from Lemma 
3.2(b) with $P,P',Z$ replaced by $P^{n-1},P'{}^{n-1},g\i(P'{}^{n-1})^{P^{n-1}}g$.

We show that $(d_n)$ holds, using $(a_n),(b_n)$. This follows from Lemma 3.2(d) with 
$P,P',Z$ replaced by $P^{n-1},P'{}^{n-1},g\i(P'{}^{n-1})^{P^{n-1}}g$.

We show that $(e_n)$ holds, using $(a_n)$. Since $\po(P'{}^n,P^n)=w_n$, $P'{}^n\sub P'$
and $w_n\in{}^JW$ we have $\po(P',P^n)=w_n$. We also have $\po(P',P_n)=u_0u_1\do u_n$ 
(see Lemma 9.3). Since $P^n=P_n$, we have $w_n=u_0u_1\do u_n$. 

This completes the inductive proof.

\proclaim{Proposition 9.5} Let $J\sub I$. Let $\ss=(J_n,J'_n,u_n)_{n\ge 0}\in S(J,\d)$
and let $\tt\in\ct(J,\d)$ be the corresponding element under the bijection in 2.4. Then
$Z_{J,\d}^\ss={}^\tt Z_{J,\d}$. 
\endproclaim
Let $P\in\cp_J,P'\in\cp_{\d(J)},g\in A_1(P,P')=A'_1(P,P')$. Assume that
$$(P,P',U_{P'}gU_P)\in Z_{J,\d}^{\ti\ss},(P,P',U_{P'}gU_P)\in{}^\tt Z_{J,\d}$$
where $\ti\ss=(\tJ_n,\tJ'_n,\tu_n)_{n\ge 0}\in S(J,\d)$. Using Lemma 9.4(e), we see 
that 

$\tu_0\tu_1\do\tu_n=u_0u_1\do u_n$
\nl
for all $n$. Using Lemma 2.5, we have $\ti\ss=\ss$. Thus, 
${}^\tt Z_{J,\d}\sub Z_{J,\d}^\ss$. Conversely, let $(Q,Q',\g)\in Z_{J,\d}^\ss$. We 
have $(Q,Q',\g)\in{}^{\tt'}Z_{J,\d}$ for a unique $\tt'\in\ct(J,\d)$. By the first part
of the proof we have $(Q,Q',\g)\in Z_{J,\d}^{\ss'}$ where $\ss'\in S(J,\d)$ corresponds
to $\tt'$ under the bijection in 2.4. We have 
$(Q,Q',\g)\in Z_{J,\d}^\ss\cap Z_{J,\d}^{\ss'}$ and the sets $Z_{J,\d}^\ss$,
$Z_{J,\d}^{\ss'}$ are either disjoint or coincide. Thus, $\ss=\ss'$ hence $\tt=\tt'$. 
We see that $Z_{J,\d}^\ss\sub{}^\tt Z_{J,\d}$. The proposition is proved.

\head 10. Example\endhead
\subhead 10.1\endsubhead
We consider an example. Let $V$ be a finite dimensional $\kk$-vector space. Let
$G=\hG=G^1=GL(V)$. Consider two $n$-step filtrations $V_*,V'_*$:
$$0=V_0\sub V_1\sub V_2\sub\do\sub V_n=V,\qua 0\sub V'_1\sub V'_2\sub\do\sub V'_n=V$$
of $V$. The {\it type} of $V_*$ is by definition the set
$$J=\{i\in[1,n-1];\dim V_k\ne i\qua\frl k\in[0,n]\}.$$
To $V_*,V'_*$ we associate two $n^2$-step filtrations $X_*,X'_*$:
$$\align&
0=X_{10}\sub X_{11}\sub X_{12}\sub\do\sub X_{1n}=X_{20}\sub X_{21}\sub X_{22}\sub 
X_{2n}=X_{30}\\&\sub\do\sub X_{n,n-1}\sub X_{nn}=V,\endalign$$
$$\align&
0=X'_{10}\sub X'_{11}\sub X'_{12}\sub\do\sub X'_{1n}=X'_{20}\sub X'_{21}\sub X'_{22}
\sub X'_{2n}=X'_{30}\\&\sub\do\sub X'_{n,n-1}\sub X'_{nn}=V,\endalign$$
where
$$X_{ij}=V_{i-1}+(V_i\cap V'_j), (i,j)\in[1,n]\T[0,n],$$
$$X'_{ij}=V'_{i-1}+(V'_i\cap V_j), (i,j)\in[1,n]\T[0,n].$$
Here the indexing set is $[1,n]\T[0,n]$ with the identifications $1n=20,\do,(n-1)n=n0$.
We have $X_{i0}=V_{i-1},X'_{i0}=V'_{i-1}$ for $i\in[1,n]$. Hence $X_*$ (resp. $X'_*$) 
is a refinement of $V_*$ (resp. $V'_*$). If the stabilizer of $V_*$ (resp. $V'_*$) is 
the parabolic $P$ (resp. $P'$) then the stabilizer of $X_*$ (resp. $X'_*$) in $G$ is 
the parabolic $P^{P'}$ (resp. $P'{}^P$). By Zassenhaus' lemma, we have a canonical 
isomorphism
$$t:X'_{ij}/X'_{i,j-1}@>\si>>X_{ji}/X_{j,i-1}\text{ for all }(i,j)\in[1,n]\T[1,n].$$
Assume that we are given a permutation $\s:[1,n]@>>>[1,n]$ and vector space
isomorphisms $\a_i:V_i/V_{i-1}@>>>V'_{\s(i)}/V'_{\s(i)-1}$ for $i\in [1,n]$. Define a 
third $n^2$-step filtration $Y_*$ (refining $V_*$):
$$\align&0=Y_{10}\sub Y_{11}\sub Y_{12}\sub\do\sub Y_{1n}=Y_{20}\sub Y_{21}\sub Y_{22}
\sub Y_{2n}=Y_{30}\\&\sub\do\sub Y_{n,n-1}\sub Y_{nn}=V,\endalign$$ 

$Y_{i0}=V_{i-1}$ for $i\in[1,n]$,

$Y_{ij}$ is the subspace of $V_i$ containing $V_{i-1}$ such that $a_i$ carries the 
subspace $Y_{ij}/Y_{i0}$ of $V_i/V_{i-1}$ onto the subspace $X'_{\s(i),j}/X'_{\s(i),0}$
of $V'_{\s(i)}/V'_{\s(i)-1}$.
\nl
The composition 
$$Y_{ij}/Y_{i,j-1}@>\a_i>>X'_{\s(i),j}/X'_{\s(i),j-1}@>t>>X_{j,\s(i)}/X_{j,\s(i)-1}$$
is an isomorphism $b_{ij}$ for $(i,j)\in[1,n]\T[1,n]$. Define a permutation 
$\t:[1,n]\T[1,n]@>>>[1,n]\T[1,n]$ by $\t(i,j)=(j,\s(i))$.

Let $\Si$ be the set of all quadruples $(V_*,V'_*,\s,a_i)$ with $V_*,V'_*$ as above (of
prescribed types) and $\s,a_i$ are as above. Then $\Si$ may be identified with a set 
$Z_{J,y,1}$ attached to $G=GL(V)$. Here $J,\Ad(y)J$ are the types of $V_*,V'_*$. Let 
$(V_*,V'_*,\s,a_i)\in\Si$. We define a sequence $(V_*^m,V'_*{}^m,\s^m,a_i^m)_{m\ge 0}$ 
of quadruples of the same kind as $(V_*,V'_*,\s,a_i)$ as follows. Set 
$(V_*^0,V'_*{}^0,\s^0,a_i^0)=(V_*,V'_*,\s,a_i)$. Assume that $m\ge 1$ and that 
$(V_*^{m-1},V'_*{}^{m-1},\s^{m-1},a_i^{m-1})$ is already defined. Then 
$(V_*^m,V'_*{}^m,\s^m,a_i^m)$ is attached to 
$(V_*^{m-1},V'_*{}^{m-1},\s^{m-1},a_i^{m-1})$ in the same way as $(Y_*,X_*,\t,b_{ij})$ 
was attached to $(V_*,V'_*,\s,a_i)$. Then $V_*^m,V'_*{}^m$ are $n^{2^m}$-step 
filtrations of $V$ and, for $m>0$, $V_*^m,V'_*{}^m$ are refinements of $V_*^{m-1}$. Let
$J_m$ (resp. $J'_m$) be the type of $V_*^m$ (resp. $V'_*{}^m$). The set of all 
$(V_*,V'_*,\s,a_i)\in\Si$ such that $J_m,J'_m$ and the relative position of 
$V_*^m,V'_*{}^m$ is specified for each $m\ge 0$ is a locally closed subvariety of 
$\Si$. Thus we obtain a partition of $\Si$ which coincides with the partition 
$Z_{J,y,1}=\sqc_\ss Z_{J,y,1}^\ss$.

\subhead 10.2\endsubhead
Let $V$ be a $\kk$-vector space of dimension $d\ge 2$. Let $\Si$ be the set of all 
quadruples $(V_1,V'_1,a,b)$ where $V_1,V'_1$ are lines in $V$, and
$a:V_1@>\si>>V'_1$, $b:V/V_1@>\si>>V/V'_1$ are isomorphisms. (This is a special case of
the situation in 10.1 where we omit the $0$ and $d$ dimensional members of a $2$-step 
filtration.) We describe explicitly in this case the partition of $\Si$ given in 10.1.
For any $k\in[1,d]$ let $\tSi_k$ be the set of all quadruples

$(0=V_0\sub V_1\sub V_2\sub\do\sub V_k,V'_1,a,b)$
\nl
where $0=V_0\sub V_1\sub V_2\sub\do\sub V_k$ is a (partial) flag in $V$, 
$\dim V_j=j$ for all $j$, $V'_1$ is a line in $V$, $a:V_1@>\si>>V'_1$,
$b:V/V_1@>\si>>V/V'_1$ are isomorphisms and

$V'_1\cap V_{k-1}=0,V'_1\sub V_k$,

$b(V_j/V_1)=(V_{j-1}+V'_1)/V'_1$ for $j\in[1,k]$.
\nl
Define $\p_k:\tSi_k@>>>\Si$ by 
$\p_k(0=V_0\sub V_1\sub V_2\sub\do\sub V_k,V'_1,a,b)=(V_1,V'_1,a,b)$. Then $\p_k$ is 
injective and $(\p_k(\Si_k))_{k\in[1,d]}$ is a partition of $\Si$ into $d$ locally 
closed subvarieties. This is a special case of the partition of $\Si$ in 10.1 and of 
the partition of $Z_{J,y,\d}$ in 8.18.

\subhead 10.3\endsubhead
Let $V$ be a $\kk$-vector space of dimension $d\ge 2$. Let $\Si$ be the set of all 
quadruples $(V_1,H,a,b)$ where $V_1$ is a line in $V$, $H$ is a hyperplane in $V$, and
$a:V_1@>\si>>V/H$, $b:V/V_1@>\si>>H$ are isomorphisms. (This is a special case of the 
situation in 10.1 where we omit the $0$ and $d$ dimensional members of a $2$-step 
filtration.) We describe explicitly in this case the partition of $\Si$ given in 10.1.
For any $k\in[1,d]$ let $\tSi_k$ be the set of all quadruples
$$(0=V_0\sub V_1\sub V_2\sub\do\sub V_k,H,a,b)$$
where $0=V_0\sub V_1\sub V_2\sub\do\sub V_{k-1}\sub V_k$ is a (partial) flag in $V$, 
$\dim V_j=j$ for all $j$, $H$ is a hyperplane in $V$, $a:V_1@>\si>>V/H$, 
$b:V/V_1@>\si>>H$ are isomorphisms and 

$V_{k-1}=V_k\cap H$,

$b(V_j/V_1)=V_{j-1}$ for $j\in[1,k]$.
\nl
Define $\p_k:\tSi_k@>>>\Si$ by 
$\p_k(0=V_0\sub V_1\sub V_2\sub\do\sub V_k,H,a,b)=(V_1,H,a,b)$. Then $\p_k$ is 
injective and $(\p_k(\Si_k))_{k\in[1,d]}$ is a partition of $\Si$ into $d$ locally 
closed subvarieties. This is a special case of the partition of $\Si$ in 10.1 and of 
the partition of $Z_{J,y,\d}$ in 8.18.

\head 11. Parabolic character sheaves on $Z_{J,y,\d}$\endhead
\subhead 11.1\endsubhead
Assume that we are in the setup of 8.18. Let $\xx=(x_1,x_2,\do,x_r)$ be a sequence in 
$W$ such that 
$$r\ge 2,x_r=y.\tag a$$
As in 4.2 we define
$$Y_\xx=\{(B_0,B_1,B_2,\do,B_r,g)\in\cb^{r+1}\T G^1;\po(B_{i-1},B_i)=x_i,i\in[1,r],
B_r={}^gB_0\}.$$

Let ${}^yY'_\xx$ be the set of all $(B_0,B_1,B_2,\do,B_{r-1},\g)$ where 
$(B_0,B_1,B_2,\do,B_{r-1})\in\cb^r$ satisfies $\po(B_{i-1},B_i)=x_i,i\in[1,r-1]$ and 
$\g\in U_{P'}\bsl A_y(P,P')/U_P$ (with $P\in\cp_J,P'\in\cp_{J'}$ given by 
$B_0\sub P,B_{r-1}\sub P'$) satisfies $\po(B_{r-1},{}^gB_0)=x_r$ for some/any $g\in\g$.
(This definition is correct since $U_P\sub B_0,U_{P'}\sub B_{r-1}$.) We have an affine 
space bundle
$$\o:Y_\xx@>>>{}^yY'_\xx,\qua
(B_0,B_1,B_2,\do,B_r,g)\m(B_0,B_1,B_2,\do,B_{r-1},U_{P'}gU_P)$$
(with $P,P'$ as above). Define
$$\Pi_\xx:{}^yY'_\xx@>>>Z_{J,y,\d},\qua\Pi_\xx(B_0,B_1,\do,B_{r-1},\g)=(P,P',\g)$$ 
where $P\in\cp_J,P'\in\cp_{J'}$ are given by $B_0\sub P,B_{r-1}\sub P'$. 

If $\cl\in\cs(T)$ (see 4.1) is such that $x_1x_2\do x_r\in W^1_\cl$ (see 4.1), the 
local system $\tcl$ on $Y_\xx$ is defined (see 4.2); it is $\o^*$ of a well defined
local system on ${}^yY'_\xx$ denoted again by $\tcl$. We set 
$$\ck_\xx^\cl=(\Pi_\xx)_!\tcl\in\cd(Z_{J,y,\d}).$$

Now assume in addition that
$$x_i\in I\sqc\{1\}\text{ for }i\in[1,r-1].\tag b$$
Let $Y'{}^\dag_\xx$ be the set of all $(B_0,B_1,B_2,\do,B_{r-1},\g)$ where 
$(B_0,B_1,B_2,\do,B_{r-1})\in\cb^r$ satisfies $\po(B_{i-1},B_i)=\{x_i,1\},i\in[1,r-1]$ 
and $\g\in U_{P'}\bsl A_y(P,P')/U_P$ (with $P\in\cp_J,P'\in\cp_{J'}$ given by 
$B_0\sub P,B_{r-1}\sub P'$) satisfies $\po(B_{r-1},{}^gB_0)=y$ for some/any $g\in\g$. 
Define
$$\Pi^\dag_\xx:Y'{}^\dag_\xx@>>>Z_{J,y,\d}\qua
\Pi^\dag_\xx(B_0,B_1,\do,B_{r-1},\g)=(P,P',\g)$$ 
where $P\in\cp_J,P'\in\cp_{J'}$ are given by $B_0\sub P,B_{r-1}\sub P'$. Now 
${}^yY'_\xx$ is an open dense subset subset of $Y'{}^\dag_\xx$ and it carries the local
system $\tcl$. The intersection cohomology complex $IC(Y'{}^\dag_\xx,\tcl)$ is a 
constructible sheaf $\bar\cl$ on $Y'{}^\dag_\xx$ (compare 4.3); we set
$$\bar\ck_\xx^\cl=(\Pi^\dag_\xx)_!\bar\cl\in\cd(Z_{J,y,\d}).$$
Now $\Pi^\dag_\xx$ is a proper morphism. (Indeed, the condition 
$\po(B_{r-1},{}^gB_0)=y$ in the definition of $Y'{}^\dag_\xx$ can be replaced by the 
closed condition $\po(B_{r-1},{}^gB_0)\le y$ since ${}^gB_0\sub{}^gP,B_{r-1}\sub P'$ 
and $\po(P',{}^gP)=y$.) Hence we may apply the decomposition theorem \cite{\BBD} and we
see that 

$\bar\ck_\xx^\cl$ {\it is a semisimple complex on $Z_{J,y,\d}.$ }

\proclaim{Proposition 11.2} Let $\cl\in\cs(T)$ and let $A$ be a simple perverse sheaf 
on $Z_{J,\d}$. The following conditions on $A$ are equivalent:

(i) $A\dsv\ck^\cl_\xx$ for some $\xx$ as in 11.1(a) with $x_1x_2\do x_r\in W^1_\cl$;

(ii) $A\dsv\ck^\cl_{x,y}$ for some $x\in W$ such that $xy\in W^1_\cl$;

(iii) $A\dsv\bar\ck^\cl_\xx$ for some $\xx$ as in 11.1(b) with 
$x_1x_2\do x_r\in W^1_\cl$.
\endproclaim
(Compare 4.4.)

\subhead 11.3\endsubhead
Let $\cc^\cl_{J,y,\d}$ be the set of isomorphism classes of simple perverse sheaves on
$Z_{J,y,\d}$ which satisfy the equivalent conditions 11.2(i)-(iii) with respect to 
$\cl$. The simple perverse sheaves on $Z_{J,y,\d}$ which belong to $\cc^\cl_{J,y,\d}$ 
for some $\cl\in\ct$ are called {\it parabolic character sheaves}; they (or their
isomorphism classes) form a set $\cc_{J,y,\d}$.

\subhead 11.4\endsubhead
Let $(P,P',\g)\in Z_{J,y,\d}$. Let $\g'=\g\cap A'_y(P,P')$. Let $z=\po(P',P)$. Let 
$(P_1,P'_1)=\a(P,P',g)$ where $g\in\g'$ (see 8.10, 8.11). We have 
$\po(P'_1,{}^gP_1)=z\i y$ and $\po(P'{}^P,P^{P'})=z$ (see 8.10). Let $w\in W$. We can 
write uniquely $w=ab,a\in W^{J'},b\in W_{J'}$. Let 
$$\cx=\{(B,B')\in\cb\T\cb;\po(B,B')=w,\po(B',{}^gB)=y,B\sub P,B'\sub P'\}.$$
Here $g\in\g'$; the choice of $g$ is irrelevant since $U_P\sub B,U_{P'}\sub B'$. Set 
$b'=y\i by\in W_{\d(J)}$. (Recall that $yW_{\d(J)}=W_{J'}y$.) Let 
$$\align\cx'=&\{(\tB,B,\tB')\in\cb^3;\po(\tB,B)=\d\i(b'),\po(B,\tB')=az,\\&
\po(\tB',{}^g\tB)=z\i y,\tB\sub P_1,\tB'\sub P'_1\}.\endalign$$
(We have automatically $B\sub P$.) Here $g\in\g'$; the choice of $g$ is irrelevant
since $U_{P'_1}\sub\tB',U_{P_1}\sub\tB$ (another choice of $g$ is in $U_{P'_1}gU_{P_1}$
by 8.11). Define 
$$\cx@>\mu>>\cx',\qua(B,B')\m(\tB,B,\tB')$$
as follows. Define $R$ by $\po(B,R)=a,\po(R,B')=b$ (we have $l(w)=l(a)+l(b)$). Define 
$\tB'$ by $\po(B,\tB')=az,\po(\tB',R)=z\i$ (we have $l(az)+l(z\i)=l(a)$). Define $S_g$ 
(in terms of $g\in\g'$) by $\po(R,S_g)=y,\po(S_g,{}^gB)=b'$; we use 
$l(b)+l(y)=l(by)=l(yb')=l(y)+l(b')$. Set $\tB={}^{g\i}S_g$. We will show below that 

(i) if $(B,B')\in\cx$ then $(\tB,B,\tB')\in\cx'$; 

(ii) $\tB$ is independent of the choice of $g$ in $\g'$.
\nl   
Assume that (i) is already established (for fixed $g\in\g'$). Let us replace $g$ by 
$g_1\in\g'$. Then $g_1=u'gu$ where $u\in U_P,u'\in U_{P'}$. We have 
$\po(R,S_{g_1})=y,\po(S_{g_1},{}^{u'gu}B)=b'$. Hence
$$\po({}^{u'{}\i}R,{}^{u'{}\i}S_{g_1})=y,\po({}^{u'{}\i}S_{g_1},{}^gB)=b'.$$
Since $\po(R,B')\in W_{J'}$ we have $R\sub P'$ hence $u'\in R$ and ${}^{u'{}\i}R=R$ and
$\po(R,{}^{u'{}\i}S_{g_1})=y,\po({}^{u'{}\i}S_{g_1},{}^gB)=b'$.
It follows that $S_g={}^{u'{}\i}S_{g_1}$. Hence
$${}^{u\i g\i u'{}\i}S_{g_1}={}^{u\i g\i u'{}\i}({}^{u'}S_g)={}^{u\i g\i}S_g
                            ={}^{u\i}\tB=\tB$$
(since by (i), we have $\tB\sub P_1\sub P$ hence $u\in\tB$). Thus (ii) is verified.

We now verify (i) (with fixed $g\in\g'$). Since $P^{P'},P'{}^P,{}^gP_1$ have a common 
Levi $L_0$, there is a canonical trijection ${}^0B\lra {}^1B\lra {}^2B$ between the 
sets of Borels of $P^{P'},P'{}^P,{}^gP_1$ respectively, defined by 
${}^0B=U_{P^{P'}}\b,{}^1B=U_{P'{}^P}\b,{}^2B=U_{{}^gP_1}\b$ where $\b$ is a Borel of 
$L_0$ or equivalently by
$$\po({}^0B,{}^1B)=z\i,\po({}^1B,{}^2B)=y,\po({}^0B,{}^2B)=z\i y$$
(any two of these three conditions implies the third).

Assume that $B,B',R,S,tB,\tB'$ are as above. We have $R\sub P'{}^P$ (since 
$a=pos(B,R)=pos(B,P')$) and $\tB'\sub P^{P'}$ (since $z\i=\po(\tB',R)=\po(P,P')$). 
Hence $\tB'\lra R\lra Y$ under the trijection above where $Y$ is a Borel of ${}^gP_1$.
(We have $\po(\tB',R)=z\i,\po(R,Y)=y,\po(\tB',Y)=z\i y$.)

Since $P',{}gP$ have a common Levi, there is a canonical bijection $B^0\lra B^1$ 
between the sets of Borels of $P',{}^gP$ respectively, defined by $\po(B^0,B^1)=y$.
Since $\po(R,Y)=y$, we have $R\lra Y$ under this bijection; since $\po(R,S)=y$, we have
$R\lra S$ under this bijection; hence $Y=S$. We see that $S\sub{}^gP_1$ and 
$\tB\sub P_1$. We have $\po(S,{}^gB)=b'$ hence $\po(\tB,B)=\d\i(b')$. 
Define 
$$\cx'@>\mu'>>\cx,\qua(\tB,B,\tB')\m(B,B')$$
as follows. Choose $g\in\g'$. Define $R,S$ by the condition that $\tB'\lra R\lra S$ 
under the canonical trijection above. Define $B'$ by the condition 
$\po(R,B')=b,\po(B',{}^gB)=y$ (we use $l(b)+l(y)=l(by)=l(yb')=l(y)+l(b')$). We will 
show below that 

(iii) if $(\tB,B,\tB')\in\cx'$ then $(B,B')\in\cx$;

(iv) $B'$ is independent of the choice of $g$ in $\g'$.
\nl
Now (iv) follows from (iii) in the same way that (ii) follows from (i). We now verify 
(iii) (with fixed $g\in\g'$). Assume that $\tB,B,\tB',R,B'$ are as above. We have 
$R\sub P'{}^P$ hence $R\sub P'$. Since $\po(R,B')\in W_{J'}$ we have $B'\sub P'$. We 
have $\po(B,\tB')=az,\po(\tB',R)=z\i,l(az)+l(z\i)=l(a)$ hence $\po(B,R)=a$. This,
together with $\po(R,B')=b,l(ab)=l(a)+l(b)$ implies $\po(B,B')=ab=w$. Thus
$(B,B')\in\cx$ and $\mu'$ is well defined.

From the definitions we see that $\mu'$ is the inverse of $\mu$.

\subhead 11.5\endsubhead
For $z\in{}^{J'}W^J$ let $Z_{J,y,\d;z}=\{(P,P',\g)\in Z_{J,y,\d};\po(P',P)=z\}$. As in
8.10, 8.11 we have a well defined map 
$$Z_{J,y,\d;z}@>>>Z_{J_1,z\i y,\d},\qua(P,P',\g)\m(P^1,P'{}^1,\g_1)$$
where $\g_1$ is given by $\g_1\sub\g\cap A'_y(P,P')$. Let ${}^yY'_{w,y}{}^z$ be the 
inverse image of $Z{J,y,\d;z}$ under the canonical map ${}^yY'_{w,y}@>>>Z_{J,y,\d}$. 
For $a,b'$ as in 11.4, let
$$Y''_{\d\i(b'),az,z\i y}
={}^{z\i y}Y'_{\d\i(b'),az,z\i y}\T_{Z_{J_1,z\i y,\d}}Z_{J,y,\d;z}$$
where the fibre product is formed using the canonical maps
$${}^{z\i y}Y'_{\d\i(b'),az,z\i y}@>>>Z_{J_1,z\i y,\d}@<<<Z_{J,y,\d;z}.$$
The results in 11.4 provide an isomorphism
$${}^yY'_{w,y}{}^z@>\si>>Y''_{\d\i(b'),az,z\i y}$$
compatible with the natural maps of the two sides into $Z_{J,y,\d;z}$. Hence in the 
cartesian diagram
$$\CD Y''_{\d\i(b'),az,z\i y}@>>>{}^{z\i y}Y'_{\d\i(b'),az,z\i y}\\
@VVV   @VVV\\
Z_{J,y,\d;z}@>>>Z_{J_1,z\i y,\d} \endCD$$
we may substitute $Y''_{\d\i(b'),az,z\i y}$ by ${}^yY'_{w,y}{}^z$ and we obtain a
cartesian diagram
$$\CD {}^yY'_{w,y}{}^z@>>>{}^{z\i y}Y'_{\d\i(b'),az,z\i y}\\
@VVV   @VVV\\
Z_{J,y,\d;z}@>>>Z_{J_1,z\i y,\d}  \endCD$$

\subhead 11.6\endsubhead
In the setup of 11.5, let $\cl,\cl'\in\cs(T)$ be such that
$$\cl'=\Ad(\d\i(b'{}\i))^*\cl$$
(we have $\d\i(b'{}\i)\in W_J$) and
$$wy\in W^1_\cl,\qua \d\i(b')az(z\i y)\in W^1_{\cl'}.$$
(These two conditions are equivalent. In general, for $v\in W$ we have 
$W^1_{\Ad(v)^*\cl}=v\i W^1_\cl\d(v)$. In our case we have

$\d\i(b')az(z\i y)=\d\i(b')(ayb')b'{}\i=\d\i(b')wyb'{}\i$.)
\nl
Let $\tcl,\tcl'$ the local systems on ${}^yY'_{w,y},{}^{z\i y}Y'_{\d\i(b'),az,z\i y}$ 
corresponding as in 11.1 to $\cl,\cl'$. From the definitions we see that the inverse
image of $\tcl'$ under ${}^yY'_{w,y}{}^z@>>>{}^{z\i y}Y'_{\d\i(b'),az,z\i y}$ (in the
cartesian diagram above) is the same as the restriction of $\tcl$ to 
${}^yY'_{w,y}{}^z$.

\subhead 11.7\endsubhead
In the setup of 11.4, we assume in addition that $b\in W_{\Ad(y)\d(J_1)}$ that is,
$b'\in W_{\d(J_1)}$ where $J_1=J\cap\d\i\Ad(y\i z)J$. Assume that $(B,B')\in\cx$. We 
show that 
$$B\sub P_1.$$
Let $R,S,\tB,\tB'$ be as in the definition of $\mu$ in terms of some $g\in\g'$. (Thus, 
$\mu(B,B')=(\tB,B,\tB')$.) By 11.4, we have $S\sub{}^gP_1$. Since 
$\po(S,{}^gB)=b'\in W_{\d(J_1)}$ and ${}^gP_1$ has type $\d(J_1)$, we have 
${}^gB\sub{}^gP_1$ and $B\sub P_1$. 

\subhead 11.8\endsubhead
In the setup of 11.4 and 11.7, we have for $B,B',\tB,\tB'$ as in 11.7
$az=\po(B,\tB')\in{}^{J_1}W$ (since $B\sub P_1$) and $\po(\tB,B)=\d\i(b')\in W_{J_1}$
hence 

$l(\d\i(b')az)=l(\d\i(b'))+l(az)$.
\nl
It follows that $\po(\tB,\tB')=\d\i(b')az$. We see that $(\tB,B,\tB')\m(\tB,\tB')$ 
defines an isomorphism $\cx'@>\si>>\cx''$ where
$$\cx''=\{(\tB,\tB')\in\cb^2;\po(\tB,\tB')=\d\i(b')az,\po(\tB',{}^g\tB)=z\i y,
                                                      \tB\sub P_1,\tB'\sub P'_1\}.$$
Combining with the isomorphism $\mu:\cx@>\si>>\cx'$ we see that $(B,B')\m(tB,\tB')$ is
an isomorphism
$$\cx@>\si>>\cx''.\tag a$$
Now let
$$Y''_{\d\i(b')az,z\i y}
={}^{z\i y}Y'_{\d\i(b')az,z\i y}\T_{Z_{J_1,z\i y,\d}}Z_{J,y,\d;z}$$
where the fibre product is formed using the canonical maps
$${}^{z\i y}Y'_{\d\i(b')az,z\i y}@>>>Z_{J_1,z\i y,\d}@<<<Z_{J,y,\d;z}.$$
Then (a) gives rise to an isomorphism
$${}^yY'_{w,y}{}^z@>\si>>Y''_{\d\i(b')az,z\i y}$$
compatible with the natural maps of the two sides into $Z_{J,y,\d;z}$. As in 11.5 this 
gives rise to a cartesian diagram
$$\CD {}^yY'_{w,y}{}^z@>>>{}^{z\i y}Y'_{\d\i(b')az,z\i y}\\
@VVV   @VVV\\
Z_{J,y,\d;z}@>>>Z_{J_1,z\i y,\d}  \endCD$$
Now let $\cl,\cl'$ be as in 11.6. Let $\tcl,\tcl'$ the local systems on 
${}^yY'_{w,y},{}^{z\i y}Y'_{\d\i(b')az,z\i y}$ corresponding as in 11.1 to $\cl,\cl'$. 
From the definitions we see that the inverse image of $\tcl'$ under 
${}^yY'_{w,y}{}^z@>>>{}^{z\i y}Y'_{\d\i(b')az,z\i y}$ (in the cartesian diagram above) 
is the same as the restriction of $\tcl$ to ${}^yY'_{w,y}{}^z$.

\subhead 11.9\endsubhead
Let $w=ab,b'$ be as in 11.4. Let $\ss=(J_n,J'_n,u_n)_{n\ge 0}\in S(J,\Ad(y)\d)$. For 
$\xx$ as in 11.1(a) let ${}^yY'_\xx{}^\ss$ be the inverse image of $Z{J,y,\d}^\ss$ 
under the canonical map ${}^yY'_\xx@>>>Z_{J,y,\d}$. In the last (cartesian) diagram in 
11.5 (with $z=u_0$), the inverse image of $Z_{J_1,z\i y,\d}^{\ss_1}$ 

-under $Z_{J,y,\d;z}@>>>Z_{J_1,z\i y,\d}$ is $Z_{J,y,\d}^\ss$,

-under ${}^{z\i y}Y'_{\d\i(b'),az,z\i y}@>>>Z_{J_1,z\i y,\d}^{\ss_1}$ is
${}^{z\i y}Y'_{\d\i(b'),az,z\i y}{}^{\ss_1}$,

-under ${}^yY'_{w,y}{}^z@>>>Z_{J_1,z\i y,\d}$ (the two compositions in the diagram) is 
${}^yY'_{w,y}{}^\ss$.

Therefore that cartesian diagram restricts to a cartesian diagram 
$$\CD  {}^yY'_{w,y}{}^\ss@>>>{}^{z\i y}Y'_{\d\i(b'),az,z\i y}{}^{\ss_1}\\
@VVV                            @VVV\\
Z_{J,y,\d}^\ss @>>>Z_{J_1,z\i y,\d}^{\ss_1}    \endCD$$

\subhead 11.10\endsubhead
In the setup of 11.9 assume in addition that $b,b'$ are as in 11.7. As in 11.9, the 
cartesian diagram in 11.8 restricts to a cartesian diagram 
$$\CD {}^yY'_{w,y}{}^\ss@>>>{}^{z\i y}Y'_{\d\i(b')az,z\i y}{}^{\ss_1}\\
@VVV                            @VVV\\
Z_{J,y,\d}^\ss @>\vt>>Z_{J_1,z\i y,\d}^{\ss_1} \endCD$$

\subhead 11.11\endsubhead
Let $\cl\in\cs(T)$, let $\xx$ be as in 11.1(a) (with $r=3$) and let 
$\ss\in S(J,\Ad(y)\d)$. Let $\cl\in\cs(T)$ be such that $x_1x_2y\in W^1_\cl$ and let 
$\tcl$ be the corresponding local system on ${}^yY'_{x_1,x_2,y}$. Let $A'$ be a simple 
perverse sheaf on ${}^yY'_{x_1,x_2,y}{}^\ss$ such that 
$A'\dsv((\Pi_{x_1,x_2,y})_!\tcl)|_{Z_{J,y,\d}^\ss}$. We show that 

(a) {\it there exists $x_0\in W$ such that $x_0y\in W^1_\cl$ and }
$A'\dsv((\Pi_{x_0,y})_!\tcl)|_{Z_{J,y,\d}^\ss}$
\nl
(we denote the local system on ${}^yY'_{x_0,y}$ corresponding to $\cl$ again by 
$\tcl$).

The proof is similar to that of Lemma 4.8. We argue by induction on $l(x_2)$. If 
$l(x_2)=0$ then $x_2=1$, we may identify ${}^yY'_{x_1,x_2,y},{}^yY'_{x_1,y}$ and the 
result is obvious. Assume now that $l(x_2)>0$. We can find $s\in W$ such that 
$l(s)=1,l(x_2)>l(sx_2)$.

Assume first that $l(x_1s)=l(x_1)+1$. Then
$${}^yY'_{x_1,x_2,y}@>\si>>{}^yY'_{x_1s,sx_2,y},\qua
(B_0,B_1,B_2,\g)\m(B_0,B'_1,B_2,\g)$$
with $B'_1$ given by 
$$\po(B_1,B'_1)=s,\po(B'_1,B_2)=sx_2\tag b$$
is an isomorphism.
Under this isomorphism, $\tcl$ on ${}^yY'_{x_1,x_2,y}$ corresponds to the analogous 
local system on ${}^yY'_{x_1s,sx_2,y}$. We have
$A'\dsv((\Pi_{x_1s,sx_2,y})_!\tcl)|_{Z_{J,y,\d}^\ss}$. We may apply the induction 
hypothesis to $x_1s,sx_2,y$; the desired result follows.

Assume next that $l(x_1s)=l(x_1)+1$. Then we have a partition
${}^yY'_{x_1,x_2,y}={}'Y\sqc{}''Y$ where ${}'Y$ is the open subset defined by
$\po(B_0,B'_1)=x_1$ (with $B'_1$ as in (b)) and ${}''Y$ is the closed subset defined by
$\po(B_0,B'_1)=x_1s$ (with $B'_1$ as in (b)). Let ${}'\Pi:{}'Y@>>>Z_{J,y,\d},
{}''\Pi:{}''Y@>>>Z_{J,y,\d}$ be the restrictions of $\Pi_{x_1,x_2,y}$ to ${}'Y,{}''Y$.
By general principles, we have either

(c) $A'\dsv({}'\Pi_!\tcl)|_{Z_{J,y,\d}^\ss}$ or 
 
(d) $A'\dsv({}''\Pi_!\tcl)|_{Z_{J,y,\d}^\ss}$
\nl
where the restriction of $\tcl$ to ${}'Y$ or ${}''Y$ is denoted again by $\tcl$.

Assume that (d) holds. Then 
$$\io'':{}''Y@>>>{}^yY'_{x_1s,sx_2,y},\qua(B_0,B_1,B_2,\g)\m(B_0,B'_1,B_2,\g)$$
with $B'_1$ as in (b), is a line bundle and $\io''_!\tcl$ is up to shift the local 
system on ${}^yY'_{x_1s,sx_2,y}$ attached to $\cl$ (we denote it again by $\tcl$). 
Since ${}''\Pi=\Pi_{x_1s,sx_2,y}\io''$, it follows that 
$A'\dsv((\Pi_{x_1s,sx_2,y})_!\tcl)|_{Z_{J,y,\d}^\ss}$. We may apply the induction 
hypothesis to $x_1s,sx_2,y$; the desired result follows.

Assume now that (c) holds. Then
$$\io':{}'Y@>>>{}^yY'_{x_1,sx_2,y},\qua(B_0,B_1,B_2,\g)\m(B_0,B'_1,B_2,\g)$$
with $B'_1$ as in (b), makes $Y'$ into the complement of the zero section of a line
bundle over ${}^yY'_{x_1,sx_2,y}$ and we have ${}'\Pi=\Pi_{x_1s,sx_2,y}\io'$. In the 
case where 

(e) {\it the inverse image of $\cl$ under the coroot $\kk^*@>>>T$ corresponding to
$y\i x_2\i sx_2y$ is $\bbq$},
\nl
(so that $x_1sx_2y\in W^1_\cl$), $\tcl$ (on ${}'Y$) is $\io'{}^*$ of the local system 
on ${}^yY'_{x_1,sx_2,y}$ (denoted by $\tcl'$) hence we have an exact triangle 
consisting of $\io'_!\tcl,\tcl'$ and a shift of $\tcl'$. Hence 
$A'\dsv((\Pi_{x_1,sx_2,y})_!\tcl)|_{Z_{J,y,\d}^\ss}$. We may apply the induction 
hypothesis to $x_1,sx_2,y$; the desired result follows. In the case where (e) does not
hold, we have $\io'_!\tcl=0$ hence ${}'\Pi_!\tcl=0$, a contradiction. (a) is proved.

\subhead 11.12\endsubhead    
Let $\ss=(J_n,J'_n,u_n)_{n\ge 0}\in S(J,\Ad(y)\d)$. For $r\gg 0$ we have 
$J_r=J'_r=J_{r+1}=J'_{r+1}=\do$, and $u_r=u_{r+1}=\do=1$. Let $P\in\cp_{J_r}$. Let 
$L^\ss$ be a Levi of $P$. Then 
$$C^\ss=\{g\in G^1;{}^gL^\ss=L^\ss,\po(P,{}^gP)=y_r\}$$
(where $y_r=u_{r-1}\i\do u_0\i y$) is a connected component of $N_{\hG}(L^\ss)$. Let 
$X$ be a character sheaf on $C^\ss$ (the definition in 4.5 is applicable since $C^\ss$ 
is a connected component of an algebraic group with identity component $L^\ss$). We 
regard $X$ as a simple perverse sheaf on $U_P\bsl A_{y_r}(P,P)/U_P$ via the obvious 
isomorphism $C^\ss@>\si>>U_P\bsl A_{y_r}(P,P)/U_P$. Now $X$ is $P$-equivariant for the 
conjugation action of $P$. Hence there is a well defined simple perverse sheaf $X'$ on 
$G\T_P(U_P\bsl A_{y_r}(P,P)/U_P)$ (with $P$ acting on $G$ by right translation) whose
inverse image under 
$$G\T(U_P\bsl A_{y_r}(P,P)/U_P)@>>>G\T_P(U_P\bsl A_{y_r}(P,P)/U_P)$$
is a shift of the inverse image of $X$ under
$$pr_2:G\T(U_P\bsl A_{y_r}(P,P)/U_P)@>>>U_P\bsl A_{y_r}(P,P)/U_P.$$
We may regard $X'$ as a simple perverse sheaf on $Z_{J_r,y_r,\d}^{\ss_r}$ via the 
isomorphism 
$$G\T_P(U_P\bsl A_{y_r}(P,P)/U_P)@>\si>>Z_{J_r,y_r,\d}^{\ss_r},\qua
(g,\g)\m({}^gP,{}^gP,g\g g\i).$$
Now let $\vt:Z_{J,y,\d}^\ss@>>>Z_{J_r,y_r,\d}^{\ss_r}$ be a composition of maps in 
8.20(a), a smooth map with connected fibres. Then $\tX=\ti\vt(X')$ is a simple perverse
sheaf on $Z_{J,y,\d}^\ss$. Let $\hX$ be the simple perverse sheaf on $Z_{J,y,\d}$ whose
support is the closure in $Z_{J,y,\d}$ of $\supp\tX$ and whose restriction to
$Z_{J,y,\d}^\ss$ is $\tX$.

Let $\cc'_{J,y,\d}{}^\ss$ be the class of simple perverse sheaves on $Z_{J,y,\d}^\ss$ 
consisting of all $\tX$ as above. Let $\cc'_{J,y,\d}$ be the class of simple perverse 
sheaves on $Z_{J,y,\d}$ consisting of all $\hX$ as above (where $\ss$ varies). The 
isomorphism classes of objects in $\cc'_{J,y,\d}$ are in bijection with the set of     
pairs $(\ss,X)$ where $\ss\in S(J,\Ad(y)\d)$ and $X$ is a character sheaf on $C^\ss$ 
(as above). 

\proclaim{Lemma 11.13} Let $\ss=(J_n,J'_n,u_n)_{n\ge 0}\in S(J,\Ad(y)\d)$. Let 
$\cl\in\cs(T)$ and let $w\in W$ be such that $wy\in W^1_\cl$; let $\tcl$ be the 
corresponding local system on ${}^yY'_{w,y}$. Let $A'$ be a simple perverse sheaf on 
$Z_{J,y,\d}^\ss$ such that $A'\dsv(\Pi_{w,y})_!\tcl|_{Z_{J,y,\d}^\ss}$. Then 
$A'\in\cc'_{J,y,\d}{}^\ss$.
\endproclaim
More generally, we show that the lemma holds when $J,\ss$ are replaced by $J_n,\ss_n$ 
for any $n\ge 0$. First we show:

(a) if the result is true for $n=1$ then it is true for $n=0$.
\nl
Let $\cl',\tcl'$ be related to $\cl,\tcl$ as in 11.6. The restrictions of $\tcl,\tcl'$
to ${}^yY'_{w,y}{}^\ss$, ${}^{z\i y}Y'_{\d\i(b'),az,z\i y}{}^{\ss_1}$ (where $z=u_0$) 
are again denoted by $\tcl,\tcl'$. In the last cartesian diagram in 11.10, $\tcl$ is 
the inverse image of $\tcl'$ under 
$${}^yY'_{w,y}{}^\ss@>>>{}^{z\i y}Y'_{\d\i(b'),az,z\i y}{}^{\ss_1}$$
(see the last sentence in 11.8). By the change of basis theorem, the direct image with 
compact support of $\tcl$ under ${}^yY'_{w,y}{}^\ss@>>>Z_{J,y,\d}^\ss$ is $\vt^*$ 
($\vt$ as in 8.19) of the direct image with compact support of $\tcl'$ under
$${}^{z\i y}Y'_{\d\i(b')az,z\i y}{}^{\ss_1}@>>>Z_{J_1,z\i y,\d}^{\ss_1}.$$
In other words, 
$$(\Pi_{w,y})_!\tcl|_{Z_{J,y,\d}^\ss}=
\vt^*((\Pi_{\d\i(b'),az,z\i y})_!\tcl'|_{Z_{J_1,z\i y,\d}^{\ss_1}}).$$
Thus, $A'\dsv\vt^*((\Pi_{\d\i(b'),az,z\i y})_!\tcl'|_{Z_{J_1,z\i y,\d}^{\ss_1}})$. 
Since $\vt$ is an affine space bundle it follows that there exists a simple perverse 
sheaf $A''$ on $Z_{J_1,z\i y,\d}^{\ss_1}$ such that $A'=\ti\vt(A'')$ and 
$A''\dsv(\Pi_{\d\i(b'),az,z\i y})_!\tcl'|_{Z_{J_1,z\i y,\d}^{\ss_1}}$. Applying 
11.11(a) for $\d\i(b'),az,z\i y,\ss_1$ instead of $x_1,x_2,y,\ss$ we see that there 
exists $x_0\in W$ such that $x_0z\i y\in W^1_{\cl'}$ and 
$A''\dsv((\Pi_{x_0,y}\tcl')|_{Z_{J_1,z\i y,\d}^{\ss_1}}$. By our assumption we have 
$A''\in\cc'_{J_1,z\i y,\d}{}^{\ss_1}$. From the definitions we have
$\ti\vt(A'')\in\cc'_{J,y,\d}{}^\ss$. Thus, (a) holds.

Similarly, if the result holds for some $n\ge 1$ then it holds for $n-1$. (The proof is
the same as for $n=1$.) In this way we see that it suffices to prove the result for
$n\gg 0$. Thus we may assume that $J_0=J'_0=J_1=J'_1=\do=J$ and $u_0=u_1=\do=1$. Then
$W_Jy=yW_{\d(J)}$. In our case, ${}^yY'_{w,y}{}^\ss\ne\em$ hence there exist
$B_0,B_1\in\cb$ such that $\po(B_0,B_1)=w$ and $B_0,B_1$ are contained in the same
parabolic of type $J$. Thus we have $w\in W_J$. Let $P\in\cp_J$ with $T\sub P$. Let $L$
be the Levi of $P$ such that $T\sub L$. Then
$$C=\{c\in G^1;{}^cL=L,\po(P,{}^cP)=y\}$$
is a connected component of $N_{\hG}(L)$. Let $Y'$ be the set of all $(\b_0,\b_1,c)$ 
where $\b,\b'$ are Borels of $L$ such that $\po(\b_0,\b_1)=w$ (position with respect to
$L$ with Weyl group $W_J$) and $c\in C$ is such that ${}^c\b_0=\b_1$. Then $P$ acts on
$Y'$ by $p:(\b_0,\b_1,c)\m({}^l\b_0,{}^l\b_1,{}^lc)$ where $l\in L,p\in lU_P$. We have
a commutative diagram
$$\CD      G\T_PY'@>\si>>{}^yY'_{w,y}{}^\ss\\
           @VVV     @VVV\\
           G\T_PC @>\si>>Z_{J,y,\d}^\ss    \endCD$$
where the upper horizontal map is 
$$(g,\b_0,\b_1,c)\m({}^gB_0,{}^gB_1,{}^gc)$$
with $B_0=\b_0U_P,B_1=\b_1U_P$, the lower horizontal map is 
$$(g,c)\m({}^gP,{}^gP,U_Pgcg\i U_P),$$
the left vertical map is $(g,\b_0,\b_1,c)\m(g,c)$ and the right vertical map is 
$\Pi_{w,y}$. This commutative diagram shows that any composition factor of 
$\op_i{}^pH^i((\Pi_{w,y})_!\tcl)$ is of the form $X'$ (notation of 11.12) where $X$ is 
a character sheaf on $C$; hence it is in $\cc'_{J,y,\d}{}^\ss$. The lemma is proved.

\proclaim{Lemma 11.14} For $\ss\in S(J,\Ad(y)\d), A\in\cc_{J,y,\d}$, we set
$A^\ss=A|_{Z_{J,y,\d}^\ss}$. Then any composition factor of $\op_i{}^pH^i(A^\ss)$ 
belongs to $\cc'_{J,y,\d}{}^\ss$.
\endproclaim
We can find $\cl\in\cs(T)$ and $\xx=(x_1,x_2,\do,x_r)$ as in 11.1(b) such that
$x_1x_2\do x_r\in W^1_\cl$ and $A\dsv\bar\ck^\cl_\xx$ (see 11.1, 11.2). Since the 
complex $\bar\ck^\cl_\xx$ is semisimple (see 11.1) we have 
$\bar\ck^\cl_\xx\cong A[m]\op K'$ for some $K'\in\cd(Z_{J,y,\d})$ and some $m\in\zz$.
Hence $\bar\ck^\cl_\xx|_{Z_{J,y,\d}^\ss}\cong A^\ss[m]\op K'_1$ for some 
$K'_1\in\cd(Z_{J,y,\d}^\ss)$. It suffices to show that, if 
$A'\dsv\bar\ck^\cl_\xx|_{Z_{J,y,\d}^\ss}$, then $A'\in\cc'_{J,y,\d}{}^\ss$. As in
\cite{\CS, 2.11-2.16} we see that there exists $\cl\in\cs(T),w\in W$ such that
$wy\in W^1_\cl$ and $A'\dsv\ck^\cl_{w,y}|_{Z_{J,y,\d}^\ss}$. Using Lemma 11.13 we have 
$A'\in\cc'_{J,y,\d}{}^\ss$. The lemma is proved.

\proclaim{Lemma 11.15} If $A\in\cc_{J,y,\d}$ then $A\in\cc'_{J,y,\d}$.
\endproclaim
Since $Z_{J,y,\d}=\sqc_{J\sub I}Z_{J,y,\d}^\ss$, we can find $\ss\in S(J,\Ad(y)\d)$ 
such that $\supp(A)\cap Z_{J,y,\d}^\ss$ is open dense in $\supp(A)$. Then
$A^\ss=A|_{Z_{J,y,\d}^\ss}$ is a simple perverse sheaf on $Z_{J,y,\d}^\ss$ and 
$A^\ss\in\cc'_{J,y,\d}{}^\ss$ (Lemma 11.14). Now $A,A^\ss$ are related just as 
$\hX,\tX$ are related in 11.12. Hence $A\in\cc'_{J,y,\d}{}^\ss$. The lemma is proved.

\proclaim{Lemma 11.16} Let $\ss=(J_n,J'_n,u_n)_{n\ge 0}\in S(J,\Ad(y)\d)$,
$C,X,X',\tX,\hX$ be as in 11.12. For any $n\ge 0$ define a simple perverse sheaf $X'_n$
on $Z_{J_n,y_n,\d}^{\ss_n}$ by $X'_n=\ti\vt(X'_{n+1})$ where 
$\vt:Z_{J_n,y_n,\d}^{\ss_n}@>>>Z_{J_{n+1},y_{n+1},\d}^{\ss_{n+1}}$ is as in 8.20(a) for
$n\ge 0$ and $X'_n=X'$ for $n\gg 0$. Define $a_n\in W^{J'_n}$ by 
$a_n\i=u_nu_{n+1}\do u_m$ for $m\gg 0$. For any $n\ge 0$ there exists $\cl_n\in\cs(T)$ 
and $b'_n\in W_{\d(J_\iy)}$ (see 2.6) such that $a_ny_nb'_n\in W^1_{\cl_n}$ and
$X'_n\dsv(\Pi_{a_ny_nb'_ny_n\i,y_n})_!\tcl_n|_{Z_{J_n,y_n,\d}^{\ss_n}}$.
\endproclaim
Assume that the result holds for $n=1$; we show that it holds for $n=0$. By assumption 
we have $X'_1\dsv(\Pi_{ayb'_1y\i z,z\i y})_!\tcl'|_{Z_{J_1,y_1,\d}^{\ss_1}}$ where
$\tcl'=\tcl_1,z=u_0,a=a_0$. (We have $a_1y_1=ay$.) We consider the cartesian diagram in
11.10 with $w=ayb'y\i$ where $\d\i(b')a=ayb'_1y\i$. (We have 
$b'=\d(ayb'_1y\i a\i)\in W_{\d(J_\iy)}$ by 2.6.) The inverse image of $\tcl'$ under 
${}^yY'_{w,y}{}^\ss@>>>{}^{z\i y}Y'_{\d\i(b')az,z\i y}{}^{\ss_1}$ is $\tcl$ for some
$\cl\in\cs(T)$ (see 11.8). Using the change basis theorem for the cartesian diagram in
11.10 we deduce that $X'_0\dsv(\Pi_{ayb'y\i,y})_!\tcl|_{Z_{J,y,\d}^\ss}$. 

The same argument shows that, if the result holds for some $n\ge 1$ then it also holds
for $n-1$. In this way it suffices to show that the result holds for $n\gg 0$. 
Replacing $\ss,n$ by $\ss_n,0$, we may assume that $J_0=J'_0=J_1=J'_1=\do=J$,
$u_0=u_1=\do=1$ and $n=0$. Let $P,L$ be as in 11.13. We can find $w\in W_J$ such that
$X\dsv(pr_3)_!\tcl$ where $pr_3:Y'@>>>C$ is defined with $Y'$ as in 11.13 in terms of 
$w$; here $\cl\in\cs(T)$ and $\tcl$ is the corresponding local system on $Y'$. Using 
the commutative diagram in 11.13 we see that 
$X'\dsv(\Pi_{yb'y\i,y})_!\tcl|_{Z_{J,y,\d}^\ss}$ where $b'=y\i wy\in W_{\d(J)}$ and 
$\tcl$ is the local system on ${}^yY'_{w,y}{}^\ss$ corresponding to $\cl$. The lemma is
proved.

\proclaim{Lemma 11.17}Let $\ss=(J_n,J'_n,u_n)_{n\ge 0}\in S(J,\Ad(y)\d)$. Define
$a\in W^{J'}$ by $a\i=u_0u_1\do u_m$ for $m\gg 0$. Let $b'\in W_{\d(J_\iy)}$. Then the
image of $\Pi_{ayb'y\i,y}:{}^yY'_{ayb'y\i,y}@>>>Z_{J,y,\d}$ is contained in
$Z_{J,y,\d}^\ss$.
\endproclaim
Let $(B_0,B_1,\g)\in{}^yY'_{ayb'y\i,y}$. Let $(P,P',\g)=\Pi_{ayb'y\i,y}(B_0,B_1,\g)$.
We have $yb'{}\i y\i\in W_{J'}$. Hence 
$$\align&\po(P',P)=\min(W_{J'}\po(B_1,B_0)W_J)=\min(W_{J'}yb'{}\i y\i a\i W_J)\\&=
\min(W_{J'}a\i W_J)=u_0.\endalign$$
Thus, $(P,P',\g)\in Z_{J,y,\d;z}$ where $z=u_0$. Let $(P_1,P'_1,\g_1)$ be the image of 
$(P,P',\g)$ under $Z_{J,y,\d;z}@>>>Z_{J_1,z\i y,\d}$ (see 11.5). By the cartesian 
diagram in 11.8, there exists 
$$(\tB_0,\tB_1,\ti\g)\in{}^{z\i y}Y'_{\d\i(b')az,z\i y}$$
such that
$$\Pi_{\d\i(b')az,z\i y}(\tB_0,\tB_1,\ti\g)=(P_1,P'_1,\g_1).$$
We have $\d\i(b'{}\i)\in W_{J_\iy}\sub W_{J_1}$ hence
$$\align&\po(P'_1,P_1)=\min(W_{J'_1}\po(\tB_1,\tB_0)W_{J_1})=
\min(W_{J'_1}z\i a\i\d\i(b'{}\i)W_{J_1})=\\&\min(W_{J'_1}u_0\i a\i W_{J_1})=u_1.
\endalign$$
Now $\d\i(b')az=(az)(z\i y)b'_1(y\i z)$ where $b'_1=(ay)\i\d\i(b')ay\in W_{\d(J_\iy)}$ 
(see 2.6). Hence in the previous argument we may replace $B_0,B_1,\g,P,P',a,y,b'$ by
$\tB_0,\tB_1,\ti\g,P_1,P'_1,au_0,u_0\i y,b'_1$ and we see that the image 
$(P_2,P'_2,\g_2)$ of $(P_1,P'_1,\g_1)$ under $Z_{J_1,u_0\i y,\d;u_1}@>>>
Z_{J_2, u_1\i u_0\i y,\d}$ satisfies $\po(P'_2,P_2)=u_2$. Continuing this process we
find that $(P,P',\g)\in Z_{J,y,\d}^\ss$. The lemma is proved.

\proclaim{Lemma 11.18} If $A\in\cc'_{J,y,\d}$ then $A\in\cc_{J,y,\d}$.
\endproclaim
Let $\ss,\tX,\hX$ be as in the proof of Lemma 11.16. We may assume that $A=\hX$. By 
Lemma 11.16 we have $\tX\dsv(\Pi_{ayb'y\i,y})_!\tcl|_{Z_{J,y,\d}^\ss}$ where $a=a_0$ 
(see 11.16), $b'\in W_{\d(J_\iy)}$ for some $\cl\in\cs(T)$ with $ayb'\in W^1_\cl$. By
Lemma 11.17, $\Pi_{ayb'y\i,y}:{}^yY'_{ayb'y\i,y}@>>>Z_{J,y,\d}$ factors through a map
$\Pi':{}^yY'_{ayb'y\i,y}@>>>Z_{J,y,\d}^\ss$ and $\tX=\Pi'_!\tcl$. Thus there exists a
simple perverse sheaf on $Z_{J,y,\d}$ whose support is the closure in $Z_{J,y,\d}$ of
$\supp(\tX)$, whose restriction to $Z_{J,y,\d}^\ss$ is $\tX$ and which is a composition
factor of $\op_i{}^pH^i((\Pi_{ayb'y\i,y})_!\tcl)$; this is necessarily $\hX$. We see 
that $\hX\in\cc_{J,y,\d}$. The lemma is proved.

\subhead 11.19\endsubhead
For $P\in\cp_J$ let $H_P$ be the inverse image of the connected centre of $P/U_P$ under
$P@>>>P/U_P$. For $P,\tP\in\cp_J$, the groups $H_P/U_P,H_{\tP}/U_{\tP}$ are canonically
isomorphic (an element $h\in G$ that conjugates $P$ into $\tP$ induces an isomorphism 
$H_P/U_P@>\si>>H_{\tP}/U_{\tP}$ that is independent of the choice of $h$). Thus we may 
identify the groups $H_P/U_P$ (with $P\in\cp_J$) with a single torus $\D_J$ independent
of $P$. Now $\D_J$ acts (freely) on $Z_{J,y,\d}$ by $\d:(P,P',\g)\m(P,P',\g z)$ where 
$z\in H_P$ represents $\d\in\D_J$ and each piece $Z_{J,y,\d}^\ss$ is $\D_J$-stable. We 
set
$$\bZ_{J,y,\d}^\ss=\D_J\bsl Z_{J,y,\d}^\ss,
\bZ_{J,y,\d}=\D_J\bsl Z_{J,y,\d}=\sqc_\ss\bZ_{J,y,\d}^\ss.\tag a$$
The action of $G$ on $Z_{J,y,\d}$ and on $Z_{J,y,\d}^\ss$ commutes with the action of 
$\D_J$ and induces an action of $G$ on $\bZ_{J,y,\d}$ and on $\bZ_{J,y,\d}^\ss$.

For $(P,P',\g)\in Z_{J,y,\d}$ we have naturally
$$H_P/U_P\sub H_{P_1}/U_{P_1}\sub\do\sub H_{P_r}/U_{P_r}\sub\do$$
hence we may identify $\D_J$ with a (closed) subgroup of the centre of $L^\ss$ (as 
above). The bijection 8.21(a) induces a bijection
$$G\bsl\bZ_{J,y,\d}\lra\sqc_{\ss\in S(J,\Ad(y)\d)}(L^\ss/\D_J)\bsl(C^\ss/\D_J).\tag b$$
Here $C^\ss/\D_J$ is the orbit space for the free action of $\D_J$ by right translation
on $C^\ss$ (restriction of the action of the centre of $L^\ss$ by right translation)
and the action of $L^\ss/\D_J$ on $C^\ss/\D_J$ is induced by the conjugation action of 
$L^\ss$ on $C^\ss$.

\head 12. Completion\endhead
\subhead 12.1\endsubhead
Assume that $P,P'$ are two parabolics of $G$ (as in 0.1) and that 
$\Ph:P/U_P@>>>P'/U_{P'}$ is an isomorphism of algebraic groups. Let $H^\Ph_{P,P'}$ be
the set of all $(f,f')\in P\T P'$ such that $\Ph$ carries the image of $f$ in $P/U_P$
to the image of $f'$ in $P'/U_{P'}$. Clearly, $H^\Ph_{P,P'}$ is a closed connected
subgroup of $G\T G$ of dimension $\dim U_P+\dim U_{P'}+\dim(P/U_P)=\dim G$ and  
$P,P',\Ph$ can be reconstructed from $H^\Ph_{P,P'}$.

Let $\cv_G$ be the (projective) variety whose points are the $\dim(G)$-dimensional Lie
subalgebras of $\Lie(G\T G)$. We have $\Lie H^\Ph_{P,P'}\in\cv_G$.

\subhead 12.2\endsubhead
Assume now that we are in the setup of 3.1. Let $J,J',y$ be as in 8.8. For 
$(P,P',\g)\in Z_{J,y,\d}$ we set $\Ph_\g=ba\i\Ad(g):P/U_P@>\si>>P'/U_{P'}$ where
$g\in\g$ and $a,b$ are the obvious isomorphisms in
$$P/U_P@>\Ad(g)>>{}^gP/U_{{}^gP}@<a<<({}^gP\cap P')/U_{{}^gP\cap P'}@>b>>P'/U_{P'}.$$
Now $\Ph_\g$ is independent of the choice of $g$ and 
$(P,P',\g)\m\Lie H^{\Ph_\g}_{P,P'}$ is an imbedding
$$\bZ_{J,y,\d}\sub\cv_G\tag a$$
(notation of 11.19(a)).

\subhead 12.3\endsubhead
In the remainder of this section we assume that $G$ is adjoint. Recall that two 
parabolics $Q,Q'$ of $G$ are said to be opposed if their intersection is a common Levi 
of $Q,Q'$. (We then write $Q\opp Q'$.) If $B\in\cb$ and $Q$ is a parabolic in $G$, we 
write $B\lt Q$ if $B$ is opposed to some Borel of $Q$.

Let $J\sub I$. Define $J^*\sub I$ by 
$\{Q; Q\opp P\text{ for some }P\in\cp_J\}=\cp_{J^*}$. Let $y_J$ be the longest element 
in $W^{\d(J)}$. If $(P,P',g)\in\cp_J\T\cp_{\d(J)^*}\T G^1$, then $P'\opp{}^gP$ if and 
only if $g\in A_{y_J}(P,P')$ (see 3.8). Let
$$\bA_{y_J}(P,P')=H_{P'}\bsl A_{y_J}(P,P')/U_P=U_{P'}\bsl A_{y_J}(P,P')/H_P,$$
$$G^1_J
=\{(P,P',\mu);P\in\cp_J,P'\in\cp_{\d(J)^*},\mu\in\bA_{y_J}(P,P')\}=\bZ_{J,y_J,\d},$$
$$\bG^1=\sqc_{J\sub I}G^1_J.$$
(notation of 11.19(a)). We define a structure of algebraic variety on $\bG^1$. We 
identify $G^1_J$ with a subvariety of $\cv_G$ by 12.2(a). Since $G$ is adjoint we have 
$G^1_I=G^1$. By \cite{\DP}, $\bG^1$ is the closure of $G^1_I=G^1$ in $\cv_G$, so that 
$\bG^1$ is a projective variety. (In \cite{\DP} it is assumed that $G=G^1$ but the
general case can be easily reduced to this special case; in fact, $\bG^1$ is isomorphic
to the analogous variety in the case $G=G^1$.) We have a partition
$$\bG^1=\sqc_{J\sub I}\sqc_{\ss\in S(J,\Ad(y_J)\d)}\bZ_{J,y_J,\d}^\ss$$
(see 11.19(a)) refining the partition $\bG^1=\sqc_{J\sub I}\G^1_J$. Putting together 
the bijections 11.19(b) we obtain a canonical bijection
$$G\bsl\bG^1\lra
\sqc_{J\sub I}\sqc_{\ss\in S(J,\Ad(y_J)\d)}(L^\ss/\D_J)\bsl(C^\ss/\D_J).\tag a$$
Here the action of $G$ on $\bG^1$ is the extension of the conjugation action of $G$ on
$G^1$ and $(L^\ss/\D_J)\bsl(C^\ss/\D_J)$ is as in 11.19.

Let $J\sub I,\ss\in S(J,\Ad(y_J)\d)$ and let $X$ be a character sheaf on $C^\ss$ that 
is equivariant for the free action of $\D_J$ by right translation. The simple perverse 
sheaf $\tX$ on $Z_{J,y_J,\d}^\ss$ (see 11.12) is $\D_J$-equivariant hence it is a shift
of the inverse image of a well defined simple perverse sheaf $\un{\tX}$ on 
$\bZ_{J,y_J,\d}^\ss$ under the canonical map $Z_{J,y_J,\d}^\ss@>>>\bZ_{J,y_J,\d}^\ss$. 
Let $\un{\hX}$ be the simple perverse sheaf on $\bG^1$ whose support is the closure in 
$\bG^1$ of $\supp\un{\tX}$ and whose restriction to $\bZ_{J,y_J,\d}^\ss$ is $\un{\tX}$.

The character sheaves on $\bG^1$ are by definition the simple perverse sheaves on 
$\bG^1$ of the form $\un{\hX}$ with $X$ as above. The character sheaves on $\bG^1$ are 
in bijection with the set of triples $(J,\ss,X)$ where $J\sub I$, 
$\ss\in S(J,\Ad(y_J)\d)$ and $X$ is a character sheaf on $C^\ss$ that is equivariant 
for the free action of $\D_J$ by right translation.

\subhead 12.4\endsubhead
We now give another definition of a structure of algebraic variety on the set $\bG^1$
which does not use $\cv_G$.
For $B,\tB\in\cb$ and $J\sub I$, let 
$$G^{1,B,\tB}_J=\{(P,P',\mu)\in G^1_J;B\lt P',\tB\lt P,{}^g(P^{\tB})\opp P'{}^B\}$$
where $g\in\mu$ (an open subset of $G^1_J$). By the the substitution 
$(P,P',\mu)\m(B_1,B'_1,\mu)$ where $B_1=P^{\tB},B'_1=P'{}^B$, we may identify 
$G^{1,B,\tB}_J$ with the set of all triples $(B_1,B'_1,\mu)$ where $B_1,B'_1\in\cb$,
$\po(B,B'_1)=y_{J^*},\po(\tB,B_1)=y_{\d\i(J)}$, $\mu\in\bA_{y_J}(P,P')$ (with 
$P\in\cp_J,P'\in\cp_{\d(J)^*}$ defined by $B_1\sub P,B'_1\sub P'$) and 
${}^gB_1\opp B'_1$ for some/any $g\in\mu$. 

In particular, $G^{1,B,\tB}_I=\{g\in G^1;{}^g\tB\opp B\}$. Moreover, 
$G^{1,B,\tB}_\em=\{(B_1,B'_1)\in\cb\T\cb;B_1\opp\tB,B'_1\opp B\}$ (we omit the 
$\mu$-component since it is uniquely determined).

Define $\r_J:G^{1,B,\tB}_I@>>>G^{1,B,\tB}_J$ by
$$\r_J(g)=({}^{g\i}P_J,{}^g\tP_J,gH_{\tP_J}g\i U_{P_J}g)$$
where $P_J\in\cp_{\d(J)}$, $\tP_J\in\cp_{J^*}$ are given by $B\sub P_J,\tB\sub\tP_J$.
This is well defined: if $B_1,B_2\in\cb$ are opposed and $P_1\in\cp_{\d(J)}$,
$P_2\in\cp_{\d(J)^*}$ satisfy $B_1\sub P_1,B_2\sub P_2$, then 
$P_1^{B_2}\opp P_2^{B_1}$; apply this to $B_1=B,B_2={}^g\tB$.

In particular, $g\m({}^{g\i}B,{}^g\tB)$ is a morphism 
$\r_\em:G^{1,B,\tB}_I@>>>G^{1,B,\tB}_\em$. We show that 

(a) $\r_\em$ {\it is a principal $\tB/U_{\tB}$-bundle.}
\nl
We fix $(B_1,B'_1)\in G^{1,B,\tB}_\em$. Let $F=\r_\em\i(B_1,B'_1)$. Since $B_1\opp\tB$ 
and $B'_1\opp B$, we can find $g_0\in G$ such that ${}^{g_0}B_1=B,{}^{g_0}\tB=B'_1$. 
Next we can find $g_1\in G^1$ such that ${}^{g_1}B=B,{}^{g_1}B'_1=B'_1$. We have 
${}^{g_1\i g_0}B_1=B, {}^{g_1\i g_0}\tB=B'_1$. Hence $g_1\i g_0\in F$ and $F\ne\em$.

We show that $F$ is a free homogeneous $\tB/U_{\tB}$-space. Since $\tB\opp B_1$, we
have canonically $\tB/U_{\tB}=\tB\cap B_1$. It suffices to show that $F$ is a free 
homogeneous $\tB\cap B_1$-space. Now $\tB\cap B_1$ acts freely on $F$ by $t:g\m gt$. 
Let $g,g'\in F$. Since ${}^{g\i}B={}^{g'{}\i}B$ and $g'g\i\in G$, we have $g'=bg$ where
$b\in B$. Since ${}^g\tB={}^{g'}\tB$ and $g\i g'\in G$, we have $g'=g\tb$ where 
$\tb\in\tB$. We have $\tb=g\i g'={}^{g\i}b\in{}^{g\i}B=B_1$. Thus, $\tb\in\tB\cap B_1$.
This proves (a).

Similarly,

(b) $\r_J$ {\it is a principal $H_{\tP_J}/U_{\tP_J}$-bundle.} ($\tP_J$ is as above.)
\nl
Using simple roots we identify $\tB/U_{\tB}=(\kk^*)^I$ (we use again that $G$ is
adjoint). Now $(\kk^*)^I$ acts on $\kk^I$ by multiplication on each factor. This may be
regarded as an action of $\tB/U_{\tB}$ on $\kk^I$. Using the principal
$\tB/U_{\tB}$-bundle $\r_\em$ we may form the associated bundle
$$X^{B,\tB}=G^{1,B,\tB}_I\T_{\tB/U_{\tB}}\kk^I$$
(a $\kk^I$-bundle over the affine space $G^{1,B,\tB}_\em$). Now $X^{B,\tB}$ is an
affine space of dimension $\dim(G)$. We have an obvious partition 
$\kk^I=\sqc_{J\sub I}(\kk^*)^J$ and each piece is $\tB/U_{\tB}$-stable. Hence there is 
a partition
$$X^{B,\tB}=\sqc_{J\sub I}X^{B,\tB}_J$$
where $X^{B,\tB}_J=G^{1,B,\tB}_I\T_{\tB/U_{\tB}}(\kk^*)^J$. We may identify 
$X^{B,\tB}_J$ with the orbit space of $G^{1,B,\tB}_I$ by $H_{\tP_J}/U_{\tP_J}$, hence 
(using (b)) with $G^{1,B,\tB}_J$. Thus, we may identify $X^{B,\tB}$ with the subset of 
$\sqc_{J\sub I}G^1_J$ so that $X^{B,\tB}\cap G^1_J=G^{1,B,\tB}_J$. We show that for any
$J$ we have
$$G^1_J=\cup_{B,\tB\in\cb}G^{1,B,\tB}_J.\tag c$$
Let $(P,P',\mu)\in G^1_J,g\in\mu$. Since ${}^gP\opp P'$, we can find Borels 
$B_1\sub P,B'_1\sub P'$ such that ${}^gB_1\opp B'_1$. We can find $B\in\cb$ such that 
$B\lt P'$ and $P'{}^B=B'_1$. We can find $\tB\in\cb$ such that $\tB\lt P$ and 
$P^{\tB}=B_1$. Then $(P,P',\g)\in G^{1,B,\tB}_J$. This proves (c).

From (c) we deduce
$$\bG^1=\cup_{B,\tB\in\cb}X^{B,\tB}.\tag d$$
We can now define a structure of algebraic variety on $\bG^1$ by declaring that 
$X^{B,\tB}$ is an open subvariety of $\bG^1$ for any $B,\tB\in\cb$. We see that $\bG^1$
has a covering by open subsets isomorphic to the affine space of dimension $\dim(G)$. 

\subhead 12.5\endsubhead
Let $V=\{(B,B',g)\in\cb\T\cb\T G^1;B'\opp{}^gB\}$. The fibre of $pr_{12}:V@>>>\cb\T\cb$
at $(B,B')$ is just $G^{1,B',B}_I$. There is a unique free action of the torus 
$\D_\em$ (see 11.19) on $V$ whose restriction to $G^{1,B',B}_I$ is the action of 
$\D_\em=B/U_B$ appearing in 12.4(a). By 12.4(a), this makes $V$ into a principal
$\D_\em$-bundle over 
$$\{(B,B',B_1,B'_1)\in\cb^4;B_1\opp B,B'_1\opp B'\}.\tag a$$
As in 12.4, we form the associated bundle
$$X=V\T_{\D_\em}\kk^I$$
(a $\kk^I$-bundle over (a)). For $(B,B')\in\cb\T\cb$, the inclusion 
$G^{1,B',B}_I\sub V$ gives rise to an inclusion  
$G^{1,B',B}\T_{\D_\em}\kk^I\sub V\T_{\D_\em}\kk^I$ that is, $X^{B',B}\sub X$. The 
subsets $X^{B',B}$ form a partition of $X$. Let $\p:X@>>>\bG^1$ be the morphism which 
induces for any $B,B'$ the identity map $X^{B',B}@>>>X^{B',B}$. We have a partition
$$X=\sqc_{J\sub I}X_J$$
where $X_J=V_I\T_{\D_\em}(\kk^*)^J$. We may identify $X_J$ with the orbit space of $X$
by $\D_{I-J}$ hence, using 12.4(b), with the set of all quintuples 
$(B_1,B'_1,B,B',\mu)$ where 
$$(B,B',B_1,B'_1)\in\cb^4,\po(B',B'_1)=y_{J^*},\po(B,B_1)=y_{\d\i(J)},
\mu\in\bA_{y_J}(P,P')$$
(with $P\in\cp_J,P'\in\cp_{\d(J)^*}$ defined by $B_1\sub P,B'_1\sub P'$) and 
${}^gB_1\opp B'_1$ for some/any $g\in\mu$. Now $\p:X@>>>\bG^1$ restricts to the
map $X_J@>>>G^1_J$ given by $(B_1,B'_1,B,B',\mu)\m(P,P',\mu)$ (with $P,P'$ as above).

For $w\in W$ we set $V^w=\{(B,B',g)\in V;\po(B,B')=w\}$, $X^w=V^w\T_{\D_\em}\kk^I$. The
sets $X^w_J=X^w\cap X_J$ form a partition of $X$. If $(B_1,B'_1,B,B',\mu)\in X^w_J$,
then
$$\po(B_1,B)=y_{\d\i(J)}\i,\po(B,B')=w,\po(B',B'_1)=y_{J^*},\po(B'_1,{}^gB_1)=y_\em$$
for $g\in\mu$. Hence $(B_1,B,B',B'_1,{}^gB_1,g)\in Y_\xx$ where 
$\xx=(y_{\d\i(J)}\i,w,y_{J^*},y_\em)$ (see 4.2). 

Let $\cy_\xx$ be the set of all $(\b_0,\b_1,\b_2,\b_3,\b_4,\g)$ where 
$(\b_0,\b_1,\b_2,\b_3,\b_4)\in\cb^5$ and $\g\in U_{P'}\bsl A_{y_J}(P,P')/U_P$ (with 
$P\in\cp_J,P'\in\cp_{\d(J)^*}$ defined by $\b_0\sub P,\b_3\sub P'$) are such that
$(\b_0,\b_1,\b_2,\b_3,\b_4,g)\in Y_\xx$ for some/any $g\in\g$.

The obvious map $Y_\xx@>>>\cy_\xx$ is an affine space bundle. The obvious map 
$\cy_\xx@>>>X^w_J$ is a principal $\D_{\d(J)}$-bundle. Let $\cl\in\cs(T)$ be such that
$y_{\d\i(J)}\i wy_{J^*}y_\em\in W^1_\cl$ and such that the associated local system 
$\tcl$ on $Y_\xx$ is the inverse image under the composition 
$Y_\xx@>>>\cy_\xx@>>>X^w_J$ of a local system $\tcl_0$ on $X^w_J$. (In any case, $\tcl$
is the inverse image under $Y_\xx@>>>\cy_\xx$ of a well defined local system on 
$\cy_\xx$ and we require that this last local system is $\D_{\d(J)}$-equivariant.) Let 
$\KK^\cl_{w,J}$ be the direct image with compact support of $\tcl_0$ under 
$X^w_J@>>>\bG^1$ (restriction of $\p:X@>>>\bG^1$). 

We give a second definition of character sheaves on $\bG^1$ as the simple perverse
sheaves on $\bG^1$ which are composition factors of $\op_i{}^pH^i(\KK^\cl_{w,J})$ for
some $w,J,\cl$ as above. We expect that this coincides with the definition in 12.3.

\subhead 12.6\endsubhead
There is a unique simple perverse sheaf $\SS$ on $G^1$ such that:

(a) $\SS$ is a direct summand of the perverse sheaf $(pr_1)_!\bbq[\dim G]$ where
$pr_1:\{(g,B)\in G^1\T\cb;{}^gB=B\}@>>>G^1$ is the first projection (a small map);

(b) if $J\subsetneqq I,\d(J)=J$, then $\SS$ is not a direct summand of the perverse
sheaf $(pr_1)_!\bbq[\dim G]$ where $pr_1:\{(g,P)\in G^1\T\cp_J;{}^gP=P\}@>>>G^1$ is the
first projection (a small map).
\nl
Let $\tSS$ be the simple perverse sheaf on $\bG^1$ such that $\tSS|_{G^1}=\SS$. For
$x\in\bG^1$ let $G_x$ be the stabilizer of $x$ in $G$ and let $\ch^i_x(\tSS)$ be the
stalk at $x$ of the $i$-th cohomology sheaf of $\tSS$. 

We conjecture that the following three conditions on $x\in\bG^1$ are equivalent:

(c) $\ch^i_x(\tSS)\ne 0$ for some $i$;

(d) $\sum_i\dim\ch^i_x(\tSS)=1$;

(e) $G_x$ is a reductive group.
\nl
If we assume that $x\in G^1$, then the equivalence of (c),(d),(e) is known.

\subhead 12.7\endsubhead
A difficulty in proving the conjecture in 12.6 is that the small map in 12.6(a) does 
not seem to extend to a small map over all of $\bG^1$. There is one case when such an 
extension exists (partially). Assume that $G=G^1=PGL(V)$ where $V$ is a $\kk$-vector 
space of dimension $d\ge 2$. Let
$$Y=\{\t\in\End(V);\dim\ker(\t)\le 1\}/\kk^*$$
where $\kk^*$ acts by scalar multiplication. For $\t$ as above let $\bar\t$ be the 
image of $\t$ in $Y$. Let $Y_0=\{\t\in\End(V);\dim\ker(\t)=1\}/\kk^*$. We may identify 
$Y$ with an open subset of $\bG^1$ so that $Y-Y_0$ corresponds to the open stratum 
$G^1$ and $Y_0$ corresponds to a codimension $1$ stratum $G^1_{J_0}$. Let $\tY$ be the 
set of all $(\bar\t,V_1\sub V_2\sub\do\sub V_d)$ where $\bar\t\in Y$ and 
$V_1\sub V_2\sub\do\sub V_d$ are subspaces of $V$ ($\dim V_i=i$) such that 
$\t(V_i)\sub V_i$ for all $i$. Then $\tY$ is smooth and $pr_1:\tY@>>>Y$ is a small map.
Its restriction to $pr_1\i(G)$ may be identified with the small map in 12.6(a). Let 
$K_\em=(pr_1)_!\bbq[\dim G]$, a perverse sheaf on $Y$. More generally, for any 
$J\sub I$ we have a perverse sheaf $K_J$ on $Y$ defined like $K_\em$ but using flags of
type $J$ (see 8.24) in $V$ instead of complete flags. Let $\tSS'=\tSS|_Y$ that is, the
simple perverse sheaf on $Y$ such that $\tSS'|_G=\SS$. Now $\SS$ is an alternating sum 
over $J$ of the perverse sheaves $K_J|_G$; hence $\tSS'$ is an alternating sum over $J$
of the perverse sheaves $K_J$. Using this, one can compute explicitly the stalks of the
cohomology sheaves of $\tSS'$ at any $x\in Y_0$. Thus one can verify that the 
conjecture in 12.6 holds in our case for $x\in Y_0$.

\enddocument